\documentclass{eccm-ecfd}

\usepackage{hyperref}

\let\OLDthebibliography\thebibliography
\renewcommand\thebibliography[1]{
  \OLDthebibliography{#1}
  \setlength{\parskip}{0pt}
  \setlength{\itemsep}{.3pt plus 0.0ex}
}

\usepackage[dvipsnames]{xcolor}

\usepackage{mdframed}

\usepackage{tikz}
\usetikzlibrary{plotmarks}
\usetikzlibrary{shapes,arrows}
\usetikzlibrary{positioning}
\usetikzlibrary{trees,mindmap}
\usetikzlibrary{fadings,shapes.arrows,shadows}
\usetikzlibrary{decorations.pathreplacing}
\usetikzlibrary{calc}
\usepackage{pgfplots,pgflibraryplotmarks}

\usepackage{moredefs}

\usepackage{amssymb,amsmath}
\let\Emptyset\emptyset

\usepackage{stmaryrd,mathabx}
\usepackage{subfig}
\usepackage{overpic}

\usepackage{tabularx}
\usepackage{multirow,multicol,booktabs}
\usepackage{paralist}

\usepackage{subfig}
\usepackage{graphicx}

\usepackage{siunitx}
\DeclareSIUnit\poise{P}
\DeclareSIUnit\darcy{d}

\newcommand{\R}{\mathbb{R}}
\newcommand{\N}{\mathbb{N}}
\renewcommand{\P}{\mathbb{P}}

\renewcommand{\d}{\mathrm d}
\newcommand{\divergence}{\operatorname{div}}
\newcommand{\trace}{\operatorname{tr}}

\newcommand{\displacementC}{u}
\newcommand{\displacement}{\boldsymbol{\displacementC}}

\newcommand{\fluxC}{q}
\newcommand{\fluxV}{q_v}
\newcommand{\flux}{\boldsymbol{\fluxC}}
\newcommand{\pressure}{p}

\newcommand{\density}{\rho}

\newcommand{\porosity}{\phi}
\newcommand{\permeability}{\boldsymbol{K}}

\newcommand{\viscosity}{\eta}

\newcommand{\elasticityC}{c}
\newcommand{\elasticity}{\boldsymbol{C}}

\newcommand{\strainC}{\epsilon}
\newcommand{\strain}{\boldsymbol{\strainC}}
\newcommand{\strainV}{\epsilon_v}

\newcommand{\stressC}{\sigma}
\newcommand{\stress}{\boldsymbol{\stressC}}
\newcommand{\stressV}{\sigma_v}

\newcommand{\stressEC}{\stressC^{\textnormal{eff}}}
\newcommand{\stressE}{\stress^{\textnormal{eff}}}

\newcommand{\traction}{\boldsymbol{t}}

\newcommand{\gravity}{\boldsymbol{g}}

\newcommand{\Id}{\boldsymbol{1}}

\newcommand{\BiotCoefficient}{\alpha}

\newcommand{\BiotModulus}{M}

\newcommand{\Kdr}{K_{\textnormal{dr}}}

\newcommand{\II}{{\color{blue}{i}}}
\newcommand{\JJ}{{\color{blue}{j}}}
\newcommand{\KK}{{\color{blue}{k}}}
\newcommand{\LL}{{\color{blue}{l}}}

\newcommand{\ii}{{\color{Green}{\kappa}}}
\newcommand{\jj}{{\color{Green}{\iota}}}

\newcommand{\Nu}{N_{\displacement}}
\newcommand{\Nuloc}{N_{\displacement,\textnormal{loc}}}

\newcommand{\Nq}{N_{\flux}}
\newcommand{\Nqloc}{N_{\flux,\textnormal{loc}}}

\newcommand{\Np}{N_{\pressure}}
\newcommand{\Nploc}{N_{\pressure,\textnormal{loc}}}

\newcommand{\s}{{\color{red}{s}}}

\newcommand{\transpose}{T}

\renewcommand{\emptyset}{\Emptyset}

\newcommand{\doi}[1]{\href{https://doi.org/#1}{doi:#1}}
\newcommand{\urn}[1]{\href{http://nbn-resolving.de/#1}{#1}}
\newcommand{\arXiv}[1]{\href{https://arxiv.org/abs/#1}{arXiv:#1}}

\numberwithin{equation}{section}

\title{A mixed discontinuous-continuous Galerkin time discretisation
for Biot's system}

\author{UWE K\"OCHER$^{1,\dagger,\ast}$ AND MARKUS BAUSE$^{1,\ddagger}$}

\heading{Uwe K\"ocher and Markus Bause}

\address{%
$^{1}$ Helmut-Schmidt-University,
University of the Federal Armed Forces Hamburg\\
Holstenhofweg 85, 22043 Hamburg, Germany\\
\texttt{https://hsu-hh.de/mbm}\\
\begin{minipage}{.51\linewidth}
$^{\dagger}$ \texttt{koecher@hsu-hamburg.de} ($^{\ast}$ corresponding author)\\
$^{\ddagger}$ \texttt{bause@hsu-hamburg.de}
\end{minipage}
}


\abstract{%
We study higher-order space-time variational discretisations
for modeling complex processes in porous media
  that include fluid and structure interactions
which are of fundamental importance in many engineering fields
with applications in subsurface processes, battery-design and biomechanics.
For the discretisation in time we deploy discontinuous Galerkin dG($r$) and
continuous Galerkin cG($r$) discretisation families.
Moreover we introduce a new coupled dG($r$)-cG($q$) mixed time discretisation and
show numerically the stability advantages
  in the case of incompatible initial data
  under massively reduced computational costs.
For the discretisation in space we use a mixed finite element method
for the flow problem to ensure local mass conservation and
a continuous Galerkin method for the mechanics.
We consider solving sequentially the coupling of flow and mechanics with
the fixed-stress iterative approach
such that we can reuse our system solver and preconditioning technologies
for the arising block system matrices from higher-order in time discretisations.
Numerical experiments show
firstly the undeniable advantages of discontinuous Galerkin time discretisations
dG(0) and dG(1) over the continuous Galerkin time discretisation cG(1)
in the case of incompatible initial data,
secondly the advantages of the new coupled dG(1)-cG(1) in time approach
in the case of incompatible initial data
with massively reduced computational costs and better accuracy
compared to the dG(1) time discretisation and
thirdly the performance and efficiency differences of the dG(0), cG(1), dG(1)
and the new dG(1)-cG(1) fixed-stress solver approaches
for a sophisticated and physically relevant three-dimensional numerical example.
}

\begin{document}
\thispagestyle{empty}

\section{Introduction}
\label{sec:1:Introduction}

Many engineering relevant fields involve interactions between flow and mechanical
deformation throughout porous media.
Important applications in environmental and petroleum engineering include
subsurface carbon sequestration, compaction drive,
hydraulic and thermal fracturing, and oil recovery for instance.
Such problems belong to the family of thermo-hydro-mechanical (THM) problems
in a more general sense. THM problems consist of coupled heat transfer (T),
fluid flow (H) and matrix deformation (M).
THM problems usually degenerate into their hydro-mechanical (HM) effects only
for consolidation in fluid-saturated porous media processes, which are also known
as poroelasticity problems.
In reservoir engineering problems, for instance carbon sequestration problems,
the characteristic time scale for the reservoir simulation is quite large compared
to the intrinsic time scale of micro-mechanical behaviour such that
dynamic effects are neglected.
Dynamic poroelasticity models provide promising possibilities
in battery-engineering for the design of next-generation batteries
as well as
in biomechanics for simulating modern vibration therapy for osteoporosis processes
of trabeculae bones or estimating stress levels induced by tumour growth within the brain.
The preceding vivid examples show the particular importance of the ability to
simulate complex processes containing coupled flow and mechanical behaviour.
However,
  mathematical modeling,
  the development of numerical approximation schemes
  and the development of efficient solver technology of such coupled processes
remain challenging tasks for current research.

In this work, we study variational space-time numerical approximations of the
coupled flow and deformation phenomena
in its quasi-static limit for the mechanical deformation.
Variational discretisation schemes for the temporal variable got a new spirit for
their application to numerous fields of engineering relevance since sophisticated
numerical solver technology was recently established
but their origin can be traced back at least to 1969; cf. e.g. \cite{Argyris1969}.
Using variational time methods offers some appreciable advantages such as
the flexibility of the location of the degrees of freedom in time and
the straightforward construction of higher-order approximations.
In addition to that, an uniform variational method in space and time appears to
be advantageous for stability and a priori error analyses of the discrete schemes;
cf. e.g. \cite{Koecher2015b}.
In \cite{Koecher2015a,Koecher2015} we studied the application of
variational time discretisations for non-stationary diffusion problems and
for the time domain elastic wave equation. Further, in \cite{Koecher2015b} we
presented error analysis for the existence and uniqueness of the semidiscrete
and the fully discrete approximation for non-stationary diffusion problems.
In \cite{Koecher2016} we presented briefly the application of a higher-order
discontinuous in time variational scheme for Biot's system of poroelasticity.
In this work, we present details how the fully discrete systems in \cite{Koecher2016}
were established and further study the numerical properties of the proposed schemes.

In \cite{Mikelic2012} it is shown that the quasi-static Biot's system can be
obtained by taking the singular limit of the more general dynamic Biot-Allard's
system. Hereby, the Biot-Allard's system is the result of a linearisation of
dynamic fluid-structure interactions at the pore level.
The further use of an homogenisation approach allows simulations at the
macroscale level, i.e. the finite element level, with statically homogeneous
material and flow properties. The resulting system of equations consist of
an elastic wave equation for describing the dynamic behaviour of the mechanical
deformation and a diffusion equation for describing the flow phenomena.

Further, solving the coupling of flow and mechanics sequentially by using
operator splitting methods are studied in the literature
and stability as well as convergence results were established; cf.
\cite{Settari1998,Jha2007,White2008,Kim2011a,Kim2011,Mikelic2013,Koecher2016a}.

General state of the art approaches exploit low-order finite difference methods
for the time integration, such as the backward Euler scheme,
the Crank-Nicolson scheme, $\theta$-schemes or BDF-schemes.
Thus, the application of unconditionally stable higher order schemes for the time
discretisation represents an innovation over previous works.

The novelity of this work is summarised by
\begin{itemize}
\item the introduction of a new coupled dG($r$)-cG($q$) mixed time discretisation
which is stable for incompatible initial data as
  the discontinuous Galerkin in time discretisation
and having the reduced computational costs as
  the continuous Galerkin in time discretisation
for $r=q \ge 1$,

\item the detailed study and derivation of fully-discrete algebraic systems
of higher-order discontinuous and continuous Galerkin time discretisations
combined with high-order multi-physics finite elements in space for
Biot's system of poroelasticitiy
for the implementation of our efficient parallel-distributed solver suites
with monolithic and iteratively coupled solvers
\cite{Koecher2018z} respectively \cite{Koecher2018y}, and

\item sophisticated three-dimensional numerical experiments of
physically relevance to show the advantages, disadvantages and performance of
the schemes and implemenations.
\end{itemize}

The further structure of this paper is organised as follows.
Firstly we derive a mathematical model of the quasi-static Biot's system
from conservation and balance laws of the continuum problem in
Sec. \ref{sec:2:QuasiStaticBiotSystem}.
Then we present the fixed-stress iterative approach for decoupling the
flow and the mechanical subproblems for the infinite-dimensional space-time
problem in Sec. \ref{sec:3:OperatorSplitting:Continuum}.
Furthermore we study in detail the application of variational space-time
discretisations on the elliptic-mixed partial differential equation (PDE) system
in Sec. \ref{sec:4:VarSpaceTimeDisc}.
For the discretisation in time we apply members of the families of discontinuous
and continuous in time Galerkin methods and further
introduce a new coupled dG($r$)-cG($q$) mixed time discretisation
in Sec. \ref{sec:4:6:dGcGq}.
For the discretisation in space we use
a mixed finite element method (MFEM) for the flow problem and
a finite element method for the mechanical problem.
Thereby, the implementation within a distributed-parallel software of
the applied space-time discretisation allows the solution of
the presented sequential implicit splitting approach.
Finally we study the stability and convergence properties of the presented
schemes with some numerical experiments in Sec. \ref{sec:5}
and summarise this work in Sec. \ref{sec:6}.

\section{Mathematical Model: Quasi-static Biot's system}
\label{sec:2:QuasiStaticBiotSystem}

Throughout this paper we assume an isothermal single-phase flow,
small mechanical deformations, an isotropic geomechanical material and
stress-independence of flow properties;
cf. \cite{Kim2011a,Kim2011}.
The further presentation is done in a dimensionless way, i.e. precisely we let
the dimension $d \in \{2,3\}$ of the spatial domain $\Omega \subset \R^d$.
The physical model is based on poroelasticity and poroelastoplasticity theories;
cf. \cite{Biot1941,Coussy1995}.
The governing equations for coupled flow and geomechanics come from the
mass conservation law and linear-momentum balance.
Under the quasi-static assumption, the governing equation for
mechanical deformation can be expressed by the elliptic partial differential equation
\begin{equation}
\label{eq:2:1:Elasticity}
- \divergence \stress = \density_b \gravity
\end{equation}
where $\divergence(\cdot)$ is the divergence operator,
$\stress$ is the Cauchy total stress tensor,
$\density_b = \porosity \density_f + (1-\porosity)\density_s$ is the
time independent bulk density,
$\porosity$, with $0 < \porosity < 1$, is the porosity,
$\density_f$ is the fluid phase mass density,
$\density_s$ is the solid phase mass density and
$\gravity$ is some generalised volume force acting on the solid and the fluid
like gravity.
The porosity $\porosity$ is defined as the ratio of the pore volume to the
bulk volume in the deformed configuration.
Under the assumption of slightly compressible fluid, the density $\density_f$ is a
positive constant.
The solid matrix media is allowed to be heterogeneous such that
the solid phase mass density $\density_s$ and the $\porosity$ may vary within the
spatial domain. Thereby, the quantities $\density_s$ and $\porosity$ are restricted
in a way such that they are piecewise constant and uniformly bounded from below
and above as given by
\begin{displaymath}
0 < \density_s^0 \le \density_s \le \density_s^* < \infty \quad \text{and}\quad
0 < \porosity^0 \le \porosity \le \porosity^* < 1\,.
\end{displaymath}
The mechanical behaviour of the porous medium is specified by a
stress-strain relationship.
We let $\strain(\displacement)$ the linearised strain tensor applied
on the displacement $\displacement$, which is defined component-wise as
\begin{displaymath}
\strainC_{\II \JJ}(\displacement) = \frac{1}{2} \bigg(
  \frac{\partial \displacementC_{\II}}{\partial x_{\JJ}} +
  \frac{\partial \displacementC_{\JJ}}{\partial x_{\II}} \bigg)\,,
\end{displaymath}
with $1 \le \II, \JJ \le d$.
The generalised Hooke's law is given by
\begin{displaymath}
\stressE(\displacement) = \elasticity\, \strain(\displacement)\,,\quad
\stressEC_{\II \JJ} =
  \elasticityC_{\II \JJ \KK \LL}\, \strainC_{\KK \LL}(\displacement)\,,
\end{displaymath}
with $1 \le \II, \JJ, \KK, \LL \le d$,
where the summation over repeated indicies is implied, this notation is
known as Einstein summation convention, and
$\elasticity$ is a piecewise constant rank-4 elasticity tensor of
the drained soil, also known as Gassman's tensor.
For isotropic material $\elasticity$ is given by
\begin{displaymath}
\elasticityC_{\II \JJ \KK \LL} = \lambda\, \delta_{\II \JJ} \delta_{\KK \LL}
+ \mu\, ( \delta_{\II \KK} \delta_{\JJ \LL} + \delta_{\II \LL} \delta_{\JJ \KK})\,,
\end{displaymath}
where $\delta_{\II \JJ}$ denotes the Kronecker symbol,
i.e. $\delta_{\II \JJ} = 1$ if $\II=\JJ$ and zero otherwise,
and where $\lambda$ and $\mu$ are the Lam\'e parameters.
These parameters can be expressed in terms of Young's $E$ modulus and
Poisson's ratio $\nu$ as
$\lambda = E \nu((1-2\nu)(1+\nu))^{-1}$ and $\mu = E(2 (1+\nu))^{-1}$.
In addition to this, Lame's second parameter $\mu$ denotes the shear modulus $G$.
Changes in the total stress and fluid pressure are related to changes in strain
and fluid content by Biot's theory; cf. \cite{Biot1941,Biot1992,Coussy1995}.
We assume further that
$\boldsymbol E : \elasticity \boldsymbol E \ge
  \varepsilon\, | \boldsymbol E |^2 + \Kdr\, (\trace \boldsymbol E)^2$,
$\varepsilon \ge 0$,
is satisfied for any symmetric rank-2 tensor $\boldsymbol E$ and
the bulk modulus $\Kdr$ of the drained soil.
From \cite{Coussy1995}, the poroelasticity equations for
the Cauchy total stress tensor $\stress$ and
the fluid mass $m$ per bulk volume
have the following form:
\begin{displaymath}
\stress(\displacement, \pressure) := 
  \stress_0
  + \elasticity\, \strain(\displacement - \displacement_0)
  - \BiotCoefficient\, (\pressure - \pressure_0)\Id\,,
\end{displaymath}
\begin{equation}
\label{eq:2:2:massbalance}
m(\displacement, \pressure) :=
  m_0
  + \density_{f,0}\, \BiotCoefficient\, \strainV(\displacement - \displacement_0)
  + \density_{f,0}\, \BiotModulus^{-1}\, (\pressure - \pressure_0)\,,
\end{equation}
where
$\stress_0$ is the reference state of the Cauchy total stress tensor,
$\displacement_0$ is the reference state of the displacement vector field,
$\BiotCoefficient$ is the Biot's coefficient,
$\pressure$ is the fluid pressure,
$\pressure_0$ is the reference state of the fluid pressure,
$\Id$ is the rank-2 identity tensor,
$m_0$ is the reference state of the fluid mass per bulk volume,
$\density_{f,0}$ is the reference state of the fluid phase mass density,
$\strainV(\displacement) = \divergence(\displacement)$ is the volumetric strain
of the displacement $\displacement$, and
$\BiotModulus$ is the Biot's modulus.
The storage coefficient of a reservoir is defined as $c_0 := \BiotModulus^{-1}$.
\textsc{Biot} studies in his original work \cite{Biot1941}
the relation of the scalar-valued model-coupling parameters
$\BiotModulus$ and $\BiotCoefficient$ to certain physical quantities
under the assumption of an incompressible fluid that may contain air bubbles.
Nowadays, the assumption of a slightly compressible fluid is generally taken,
which is similar to the assumptions of \textsc{Biot} in the sense that the
drained soil is not completely saturated due to microscopic inclusions of air bubbles.
Following further \cite{Detournay1993,Coussy1995},
we have the micromechanical relations
\begin{equation}
\label{eq:2:3:Biot:ModulusAndCoefficient}
\BiotModulus^{-1} =
  \porosity\, c_f
  + K_s^{-1} (\BiotCoefficient - \porosity)\,,\quad \textnormal{with}\quad
c_f = K_f^{-1}\,,\quad \textnormal{and}\quad
\BiotCoefficient = 1 - K_s^{-1} \Kdr\,,
\end{equation}
where $c_f$ is the fluid compressibility,
$K_f$ is the bulk modulus of the fluid,
$K_s$ is the bulk modulus of the solid grain, and
$\Kdr$ is the drained bulk modulus.
The most recent and ongoing research
is on re-identifying the model-coupling parameters
$\BiotModulus$ and $\BiotCoefficient$ as rank-2 tensor-valued quantities
derived by homogenisation techniques and verified with experimental data.
The strain and stress tensors are decomposed in terms of their volumetric and
deviatoric parts as given by
\begin{displaymath}
\strain(\displacement) = \frac{1}{3} \strainV(\displacement) \Id
  + \boldsymbol e(\displacement)\,,\quad \textnormal{with}\quad
\strainV(\displacement) = \trace \strain(\displacement) = \divergence \displacement\,,\quad
\textnormal{and}
\end{displaymath}
\begin{displaymath}
\stress = \stressV \Id + \boldsymbol s\,,\quad \textnormal{with}\quad
\stressV = \frac{1}{3} \trace \stress\,,
\end{displaymath}
respectively, where
$\strainV(\displacement)$
is the volumetric strain, i.e. the trace of the strain tensor,
$\boldsymbol e$ is the deviatoric part of the strain tensor,
$\stressV$ is the volumetric (mean) total stress
and $\boldsymbol s$ is the deviatoric total stress tensor.
We use the convention of positive tensile and compression stresses.
The relationship between the volumetric stress and strain is given by
\begin{equation}
\label{eq:2:4:VolStressStrainRelation}
\stressV := \stressV(\displacement, \pressure) = \displaystyle
  \stressC_{v,0}
  + \Kdr\, \strainV(\displacement - \displacement_0)
  - \BiotCoefficient (\pressure - \pressure_0)\,.
\end{equation}

We further derive Biot's model of poroelasticity under quasi-static conditions
in which the stress evolution and the fluid acceleration
  due to dynamical mechanical acceleration terms $\partial_{tt} \displacement$
can be neglected since the assumptions of
$| \partial_t \nabla \pressure | \gg | \partial_{tt} \displacement |$ and
$| \partial_{tt} \displacement | = \mathcal{O}(0)$
hold for the time scale of physical interest.
We assume that the pores are small enough such that the flow through the
porous medium is laminar.
The fluid velocity $\flux$ relative to the solid phase is given by Darcy's law
\begin{equation}
\label{eq:2:5:DarcyLaw}
\flux := \displaystyle - \frac{1}{\viscosity} \permeability
  (\nabla \pressure - \density_f \gravity)\,,
\end{equation}
where $\viscosity$ is the fluid viscosity and
$\permeability$ is the permeability tensor.
The fluid viscosity is restricted to be piecewise constant and
the permeability tensor field $\permeability$ is supposed to be a piecewise
constant and being a symmetric and uniformly positive definite matrix.
For simplicity we assume that the fluid viscosity and permeability do not depend
on the history since otherwise memory terms involving fractional time-derivatives
would appear in the model; cf. \cite{Mikelic2012}.
The fluid mass conservation law under the assumption of small deformation is
given by
\begin{equation}
\label{eq:2:6:FluidMassConservation}
\partial_t m(\displacement,\pressure)
+ \nabla \cdot ( \density_{f,0}\, \flux )
= \density_{f,0}\, f\,,
\end{equation}
where $( \density_{f,0}\, \flux )$ denote the fluid mass flux,
i.e. the fluid mass flow rate per unit area and time,
and $f$ denote a volumetric source term.
More precisely, $f$ can be physically interpreted as a volumetric
fluid producing source or fluid absorbing sink term.
Additional convection terms and reaction terms are neglected for simplicity.

The fully coupled problem,
which is a parabolic partial differential equation coupled
to the quasi-static mechanical problem of Eq. \eqref{eq:2:1:Elasticity},
is derived by combining
the fluid mass conservation law given by Eq. \eqref{eq:2:6:FluidMassConservation}
with the mass balance law from Eq. \eqref{eq:2:2:massbalance}
and the relationship from Eq. \eqref{eq:2:4:VolStressStrainRelation}
as follows
\begin{displaymath}
\displaystyle
\frac{\BiotCoefficient}{\Kdr} \partial_t \stressV(\displacement, \pressure)
+ \bigg( \frac{1}{\BiotModulus}
  + \frac{\BiotCoefficient^2}{\Kdr} \bigg) \partial_t \pressure
+ \divergence \flux = f\,.
\end{displaymath}

In the sequel of this work, we use standard notation for function spaces and
their respective norms and define the function spaces $W:=L^2(\Omega)$ as
function space of square integrable functions,
$\boldsymbol V := \boldsymbol H(\operatorname{div}; \Omega) =
\{ \boldsymbol v \in L^2(\Omega)^d \,|\, \nabla \cdot \boldsymbol v \in L^2(\Omega) \}$
as function space of square integrable vector-valued functions with
square integrable divergence,
$\boldsymbol V_0 := \{ \boldsymbol v \in \boldsymbol H(\operatorname{div}; \Omega)
\,|\, (\boldsymbol v \cdot \boldsymbol n)|_{\Gamma_{\fluxC}} = 0 \}$
fulfilling additionally homogeneous boundary conditions, and
$\boldsymbol H_0 := H^1_0(\Gamma_{\displacement}; \Omega)^d$ as function space of
vector-valued componentwise $H^1(\Omega)$ functions with vanishing trace on
$\Gamma_{\displacement}$ for brevity.
Thereby, we denote by $H^m(\Omega)$ the Sobolev space of $L^2(\Omega)$ functions
with derivatives up to order $m$ in $L^2(\Omega)$.
We let $(\cdot, \cdot)$ denote the respective $L^2(\Omega)$-, $L^2(\Omega)^d$- or
$L^2(\Omega)^{d \times d}$-inner product.
Let $X$ a reflexive Banach space. We introduce the
following Bochner spaces with values in a Banach space $X$ as
\begin{displaymath}
\begin{array}{rcl}
C(\widebar I; X) &=& \{ v : [0,T] \to X\, |\, v\, \text{is continuous on $[0,T]$ with respect to $t$.} \}\,,\\[1.5ex]
L^2(I; X) &=& \bigg\{ v : (0,T) \to X \,\bigg|\, \Big( \displaystyle \int_I \| v(t) \|_X^2 \,\d t
  \Big)^{0.5} < \infty \bigg\}\,,\\[3ex]
H^1(I; X) &=& \{ v : L^2(I;X) \,|\, \partial_t v \in L^2(I;X) \}\,.
\end{array}
\end{displaymath}

\textbf{Quasi-static Biot's system.}
The fully coupled partial differential equations system in strong form,
which is precisely a coupled hydro-mechanical (HM) problem for
the approximation of the displacement, flux and pressure variables,
reads as: Find $\displacement$, $\flux$ and $\pressure$ in either $d=2$ or $d=3$
space dimensions from
\begin{equation}
\label{eq:2:7:BiotQuasiStaticSystem}
\begin{array}{rcl@{\quad}l}
\displaystyle
- \nabla \cdot \Big( \stress_0
  + \elasticity\, \strain(\displacement - \displacement_0)
  - \BiotCoefficient\, (\pressure - \pressure_0)\Id
  \Big)
  &=& \density_b\, \gravity
& \textnormal{in}\quad \Omega \times I\,,\\[3ex]
%
\displaystyle
\viscosity\, \permeability^{-1} \flux + \nabla \pressure
&=& \density_f\, \gravity
& \textnormal{in}\quad \Omega \times I\,,\\[1.5ex]
%
\displaystyle
- \frac{\BiotCoefficient}{\Kdr} \partial_t \stressV(\displacement, \pressure)
- \divergence \flux
- \bigg( \frac{1}{\BiotModulus}
  + \frac{\BiotCoefficient^2}{\Kdr} \bigg) \partial_t \pressure
&=& - f
& \textnormal{in}\quad \Omega \times I\,,\\[3ex]
\end{array}
\end{equation}
with the relationship of the volumetric stress and strain
\begin{displaymath}
\stressV(\displacement, \pressure) = \displaystyle
  \stressC_{v,0}
  + \Kdr\, \strainV(\displacement - \displacement_0)
  - \BiotCoefficient\, (\pressure - \pressure_0)\,,
\end{displaymath}
and satisfying the initial conditions
\begin{displaymath}
\begin{array}{rcl@{\quad}l @{\quad} rcl@{\quad}l}
\stress(\boldsymbol x, 0) &=& \stress_0
  & \textnormal{in}\quad \Omega \times \{0\}\,, &
\displacement(\boldsymbol x, 0) &=& \displacement_0
  & \textnormal{in}\quad \Omega \times \{0\}\,,\\[1.5ex]
\pressure(\boldsymbol x, 0) &=& \pressure_0
  & \textnormal{in}\quad \Omega \times \{0\}\,,
\end{array}
\end{displaymath}
and the boundary conditions
\begin{displaymath}
\begin{array}{rcl@{\quad}l @{\quad} rcl@{\quad}l}
\displacement &=& \displacement_D
  & \textnormal{on}\quad \Gamma_{\displacement} \times I\,, &
\stress\, \boldsymbol n &=& \traction_N
  & \textnormal{on}\quad \Gamma_{\traction} \times I\,,\\[1.5ex]
\pressure &=& \pressure_D
  & \textnormal{on}\quad \Gamma_{\pressure} \times I\,, &
\flux \cdot \boldsymbol n &=& \fluxV
  & \textnormal{on}\quad \Gamma_{\fluxC} \times I\,.
\end{array}
\end{displaymath}

The existence and uniqueness of solutions as well as regularity theory for
Biot's system of poroelasticity \eqref{eq:2:7:BiotQuasiStaticSystem} is presented
in a Hilbert space setting in \cite{Showalter2000}.
The the well-posedness is further shown in \cite{Showalter2004} for a wider
class of diffusion problems in poro-plastic media.
In \cite[Thm. 3]{Mikelic2012} the existence and uniqueness of solutions are
proofed by applying modern model homogenisation techniques which we will further
briefly summarise.
Under periodic boundary conditions for $\displacement$ and $\pressure$ with
period $L$, for smooth $L$-periodic functions $\pressure_0$, $f$ and $\stress_0$
and assumptions on the media parameters for flow and mechanics,
Biot's system from Eq. \eqref{eq:2:7:BiotQuasiStaticSystem} admits an unique
variational solution
\begin{displaymath}
\{ \displacement, \pressure \} \in
  H_{\text{per},0}^1( (0,T) \times \Omega^L  )^d \cap
  H^2(0,T; L^2(\Omega^L) )^d
\times
  H_0^1(0,T; L^2(\Omega^L)) \cap L^2(0,T; H_{\text{per}}^1(\Omega^L)) \,,
\end{displaymath}
in the poroelastic cube $\Omega^L$,
where the function spaces with the subscript ``per'' must hold periodic boundary
conditions in the spatial domain
and the function spaces with the subscript ``$0$'' must hold homogeneous
initial conditions in the temporal domain of the tensor product space-time
cube $Q = (0,T) \times \Omega^L$.
Further, additional high-order smoothness in time, i.e.
\begin{displaymath}
\{ \displacement, \pressure \} \in
  H^k(0,T; H_{\text{per}}^1(\Omega^L))^{d+1}\,,\quad
\displaystyle
\int_{\Omega^L} \displacement(\boldsymbol x,t) \,\d \boldsymbol x = 0\quad
\text{for a.e. } t \in (0,T) \,,\quad
\text{for all } k \in \N\,,
\end{displaymath}
of the quasi-static Biot's system two-field solution $\{ \displacement, \pressure \}$
was established in \cite[Lem. 3]{Mikelic2012} under the additional assumption of
$\gravity \in C^\infty_0(\R^+; L^2(\Omega^L)^d)$.
For existence and uniqueness of solutions involving more general
boundary conditions we refer to \cite{Phillips2007} and the given references therein.
Following \cite{Mikelic2012}, the quasi-static Biot's system is valid for almost
every $\boldsymbol x \in \Omega^L \subset \Omega$ in the poroelastic cube
$\Omega^L = (0,L)^d$, $d=2,3$, and for almost every $t \in I = (0,T)$, $T < \infty$.
Additionally, assume the partition of the poroelastic cube
$\Omega^L = \Omega_s \cup \Omega_f$, $\Omega_s \ne \emptyset$, $\Omega_f \ne \emptyset$,
into its deformable solid skeleton part $\Omega_s$ and
its pore part $\Omega_f$ containing the fluid on the microscale.
$\Omega_f$ is supposed to be an open periodic set with positive measure and
a smooth boundary;
cf. for details \cite{Mikelic2012,Allaire1989}.
Additionally, the pore space $\Omega_f$ is assumed to be completely filled by
a fluid at all times such that the fluid and solid phases are fully connected.
We remark that the assumptions taken on the microscale are crucial for the
development of a macroscale model by homogenisation techniques which is further
used as our simulation scale.
For well-posedness of the problem, we assume the partition of the boundary
$\partial \Omega = \Gamma_{\displacement} \cup \Gamma_{\traction}$,
with $\Gamma_{\displacement} \cap \Gamma_{\traction} = \emptyset$,
for the mechanics and the independent partition of the boundary
$\partial \Omega = \Gamma_{\pressure} \cup \Gamma_{\fluxC}$,
with $\Gamma_{\pressure} \cap \Gamma_{\fluxC} = \emptyset$ for the flow problem.
We remark that the boundary conditions for the mechanical problem may be mixed
componentwise for physically relevant problems and inhomogeneous Dirichlet
boundary values $\displacement_D$ can be incorporated in a standard way.
The initial state of strains $\strain(\displacement(\cdot, 0) - \displacement_0)$
is, by definition, homogeneous in the given model above.
The proper initialisation of the geomechanical problem is non-trivial and such,
we refer to the literature; cf. e.g. \cite{Fredrich2000}.
Being more precisely, the initial stress $\stress_0$ should satisfy
the mechanical equilibrium, i.e.
$-\nabla \cdot \stress_0 = \density_b\, \gravity(\cdot, 0)$,
and reflect the history of stress paths in the formation of the reservoir
in order to avoid stability issues.

\vskip1ex
\textbf{Weak formulation for quasi-static Biot's system.}
The fully coupled partial differential equations system in weak form reads as:
Find
$\displacement \in L^2(I; H^1_0(\Gamma_{\displacement}; \Omega)^d)
 \cap H^1(I; H^1(\Omega)^d)$,
$\flux \in L^2(I; \boldsymbol V_0)$ and
$\pressure \in H^1(I; L^2(\Omega))$ in either $d=2$ or $d=3$ space dimensions from
\begin{equation}
\label{eq:2:8:WeakForm}
\begin{array}{l}
\displaystyle
\int_I \int_\Omega \strain(\boldsymbol z) : \Big( \stress_0
  + \elasticity\, \strain(\displacement - \displacement_0)
  - \BiotCoefficient\, (\pressure - \pressure_0)\Id
  \Big) \,\d \boldsymbol x \,\d t = \\[3ex]
\displaystyle
\quad
  \int_I \int_\Omega \boldsymbol z \cdot \density_b\, \gravity \,\d \boldsymbol x \,\d t
+ \int_I \int_{\Gamma_{\traction}} \boldsymbol z \cdot \traction_N \,\d \boldsymbol x \,\d t\,,
\quad \forall \boldsymbol z \in L^2(I; H^1_0(\Gamma_{\displacement}; \Omega)^d)\,, \\[4.5ex]
%
\displaystyle
\int_I \int_\Omega \boldsymbol v \cdot \viscosity\, \permeability^{-1} \flux
  \,\d \boldsymbol x \,\d t
- \int_I \int_\Omega (\nabla \cdot \boldsymbol v)\, \pressure \,\d \boldsymbol x \,\d t
= \\[3ex]
\displaystyle
\quad
  \int_I \int_\Omega \boldsymbol v \cdot \density_f\, \gravity \,\d \boldsymbol x \,\d t
- \int_I \int_{\Gamma_{\pressure}} (\boldsymbol v \cdot \boldsymbol n)\, \pressure_D
    \,\d \boldsymbol x \,\d t\,, \quad
\forall \boldsymbol v \in L^2(I; \boldsymbol V_0)
\,,\\[4.5ex]
%
\displaystyle
- \int_I \int_\Omega w\, \divergence \flux \,\d \boldsymbol x \,\d t
- \int_I \int_\Omega w\, \bigg( \frac{1}{\BiotModulus}
  + \frac{\BiotCoefficient^2}{\Kdr} \bigg) \partial_t \pressure \,\d \boldsymbol x \,\d t
\\[3ex]
\displaystyle
\quad
- \int_I \int_\Omega w\, \frac{\BiotCoefficient}{\Kdr}
  \partial_t \stressV(\displacement, \pressure) \,\d \boldsymbol x \,\d t
= - \int_I \int_\Omega w\, f  \,\d \boldsymbol x \,\d t\,,
\,\, \forall w \in L^2(I; L^2(\Omega))\,,
\end{array}
\end{equation}
with the relationship of the volumetric stress and strain
\begin{displaymath}
\stressV(\displacement, \pressure) = \displaystyle
  \stressC_{v,0}
  + \Kdr\, \strainV(\displacement - \displacement_0)
  - \BiotCoefficient\, (\pressure - \pressure_0)\,,
\end{displaymath}
and satisfying the initial conditions
\begin{displaymath}
\begin{array}{rcl@{\quad}l @{\quad} rcl@{\quad}l}
\stress(\boldsymbol x, 0) &=& \stress_0
  & \textnormal{in}\quad \Omega \times \{0\}\,, &
\displacement(\boldsymbol x, 0) &=& \displacement_0
  & \textnormal{in}\quad \Omega \times \{0\}\,,\\[1.5ex]
\pressure(\boldsymbol x, 0) &=& \pressure_0
  & \textnormal{in}\quad \Omega \times \{0\}\,.
\end{array}
\end{displaymath}
%


\section{Operator splitting for the continuum problem}
\label{sec:3:OperatorSplitting:Continuum}

In this section we present an operator splitting method for the given continuum problem,
i.e. the quasi-static Biot's system derived in Sec. \ref{sec:2:QuasiStaticBiotSystem},
which is of the fixed-stress iterative approach type.
The goal is to decouple the mechanical problem and the flow problem.
With that, it is possible to solve the two subproblems sequentially within
an iterative scheme.
The computational overhead for the outer iteration can be resolved since
optimised and high-performance solvers can be used for each subproblem.
For low-order time discretisations the first half-step is sometimes denoted
as a half-time step but for higher order time discretisations, which have at
least two sets of unknown degrees of freedom in time per variable,
such an interpretation does not hold.
Thus, we do not use the notation of a half-time step.

The fixed-stress operator splitting method is already well studied and compared
against other splitting methods for standard time discretisations in the recent
research work of \cite{Mikelic2014,Wang2014,Mikelic2013,Kim2011}.
It is motivated by fixing the ratio of the volumetric stress
$\delta \partial_t \stressV(\displacement,\pressure) = 0$
during the first half of each fixed-point iteration step in which the flow (H)
problem is solved independently, and then the mechanical (M) problem is solved
independently with the updated values and the prescribed pressure $\pressure$ in
the second half of each iteration step.

Let `$\s$' denote the fixed-point iteration index.
The ratio of the temporal change of the volumetric stress
$\delta \partial_t \stressV(\displacement,\pressure) =
\partial_t \stressV^{\s+0.5} - \partial_t \stressV^{\s} = 0$
is constrained to vanish in the first half-step.
By using the relation from Eq. \eqref{eq:2:4:VolStressStrainRelation} we find
that
\begin{displaymath}
\partial_t \stressV^{\s+0.5} \stackrel{!}{=} \partial_t \stressV^{\s}
= \Kdr^*\, \partial_t \strainV(\displacement^{\s})
  - \BiotCoefficient\, \partial_t \pressure^{\s}
\end{displaymath}
must be fulfilled.
As a remark, we have here replaced the bulk modulus of the drained soil $\Kdr$
by some mathematically chosen $\Kdr^* > 0$.
The arbitrary choice of $\Kdr^*$ is studied in the literature for optimality in
\cite{Both2017,Koecher2016a,Mikelic2014,Mikelic2013},
but the optimal value depend on various parameters as it is studied in
\cite{Koecher2017a}.
Fluxes and pressures are prescribed
with the values $\{ \flux^{\s+0.5}, \pressure^{\s+0.5} \}$ from the first half
iteration step on the second half of each iteration step such that
$\{ \flux^{\s+1}, \pressure^{\s+1} \} := \{ \flux^{\s+0.5}, \pressure^{\s+0.5} \}$.
The fixed-stress iterative approach for the continuum problem reads as:
For $\s = 0,\,\dots,$ until convergence do:
\begin{enumerate}
\itemsep1ex

\item[(H):] {Find $\flux^{\s+1} \in L^2(I; \boldsymbol V_0)$
and $\pressure^{\s+1} \in H^1(I; L^2(\Omega))$ from
\begin{equation}
\label{eq:3:1:FixedStress:Flow}
\begin{array}{@{\!\!\!\!}l@{\;}l}
~&
\displaystyle
  \int_I (\boldsymbol v, \viscosity\, \permeability^{-1} \flux^{\s+1}) \,\d t
- \int_I (\nabla \cdot \boldsymbol v, \pressure^{\s+1}) \,\d t = \\[3.0ex]
~&
\displaystyle
\quad
  \int_I (\boldsymbol v, \density_f\, \gravity) \,\d t
- \int_I (\boldsymbol v \cdot \boldsymbol n, \pressure_D)_{\Gamma_{\pressure}} \,\d t\,,\quad
\forall \boldsymbol v \in L^2(I; \boldsymbol V_0)
\,,\\[4.5ex]
%
- &
\displaystyle
  \int_I (w, \divergence \flux^{\s+1}) \,\d t
- \int_I \bigg(w, \bigg( \frac{1}{\BiotModulus}
  + \frac{\BiotCoefficient^2}{\Kdr^*} \bigg) \partial_t \pressure^{\s+1} \bigg) \,\d t
= - \int_I (w, f) \,\d t \\[3.0ex]
~&
\displaystyle
\quad
+ \int_I (w, \BiotCoefficient\, \strainV(\partial_t \displacement^{\s})) \,\d t
- \int_I \bigg(w, \frac{\BiotCoefficient^2}{\Kdr^*}\, \partial_t \pressure^{\s} \bigg) \,\d t\,,
\quad \forall w \in L^2(I; L^2(\Omega))\,,
\end{array}
\end{equation}
satisfying the initial condition $\pressure^{\s+1}(\cdot, 0) = \pressure_0$.
}

\item[(M):] {%
Find $\displacement^{\s+1} \in L^2(I; H^1_0(\Gamma_{\displacement}; \Omega)^d)
\cap H^1(I; H^1(\Omega)^d)$ from
\begin{equation}
\label{eq:3:2:FixedStress:Mech}
\begin{array}{l}
\displaystyle
\int_I \Big( \strain(\boldsymbol z), \stress_0 \Big) \,\d t +
\int_I \Big( \strain(\boldsymbol z), \elasticity\, \strain(\displacement^{\s+1}
  - \displacement_0) \Big) \,\d t =
\int_I \Big( \strain(\boldsymbol z), \BiotCoefficient\, (\pressure^{\s+1} - \pressure_0)\Id
  - \stress_0 \Big) \,\d t\\[3ex]
\displaystyle
\quad
+ \int_I ( \boldsymbol z, \density_b\, \gravity ) \,\d t
+ \int_I ( \boldsymbol z, \traction_N )_{\Gamma_{\traction}} \,\d t\,,
\quad \forall \boldsymbol z \in L^2(I; H^1_0(\Gamma_{\displacement}; \Omega)^d)\,,
\end{array}
\end{equation}
satisfying the initial condition
$\displacement^{\s+1}(\cdot, 0) = \displacement_0$.
}
\end{enumerate}

We remark that the convergence in terms of a fixed-point iteration for the
continuum problem means here
$\{ \displacement^{\star}, \flux^{\star}, \pressure^{\star} \} =
\{ \displacement^{\s+1}, \flux^{\s+1}, \pressure^{\s+1} \} =
\{ \displacement^{\s}, \flux^{\s}, \pressure^{\s} \}$
for almost every $(\boldsymbol x,t) \in \Omega \times I$.
Also we remark that the time discretisations of (H) and (M) must not necessarily
be the same.

\section{Variational space-time discretisation}
\label{sec:4:VarSpaceTimeDisc}

In this section we present details of the applied variational space-time
discretisations. For the discretisation in time we use
a discontinuous Galerkin dG($r$) method with piecewise polynomials
  of order $r \ge 0$ in time
and a continuous Galerkin cG($q$) method with piecewise polynomials
  of order $q \ge 1$ in time.
For discontinuous Galerkin dG($r$) in time methods we use the same spaces for
trial and test functions
while we use precisely a (continuous) Petrov-Galerkin cG($q$) in time method
with continuous trial functions and discontinuous test functions;
cf. for details \cite{Koecher2015,Koecher2016,Koecher2015a}.
We are allowed to solve (efficient) time marching schemes instead of the
global space-time scheme due to the discontinuous in time test function space.

The global time interval $I=(0,T)$ is decomposed as
$\widebar I = \bigcup_{n=1}^N \widebar I_n$, $I_n = (t_{n-1}, t_n)$,
$n=1,\dots,N$, with $0 =: t_0 < \cdots < t_n < \cdots < t_N := T$.
The time subinterval length of $I_n$ is $\tau_n = t_n - t_{n-1}$ and the global
time discretisation parameter is defined as $\tau := \max_{n=1,\dots,N} \tau_n$.
The reference interval in time is defined as $\widehat I := [0,1]$.
The affine mapping $\mathcal{T}_n : \widehat I \to \widebar I_n$ and its inverse
$\mathcal{T}_n^{-1} : \widebar I_n \to \widehat I$ are given by
\begin{equation}
\label{eq:4:1}
\mathcal{T}_n( \widehat I ) := \tau_n \widehat t + t_{n-1}\,,\quad
\mathcal{T}_n^{-1}( t ) := (t - t_{n-1})/\tau_n\,,\quad
n=1,\dots,N\,.
\end{equation}
For time-discrete globally discontinuous functions, we allow the left-hand side
values $w(t_n^-)$ to be different from the right-hand side values $w(t_n^+)$.
We let $\P_r(I_n;X)$ denote the space of polynomials in time up to degree
$r \ge 0$ on $I_n$ with values in a Banach space $X$.
We introduce for some Banach space $X$
the discontinuous in time function space
\begin{equation}
\label{eq:4:2:dGr_Space}
\begin{array}{r@{\,}c@{\,}l@{\,}c@{\,}l}
\mathcal V^{\textnormal{dG($r$)}}_\tau(X) &:=&
\Big\{v \in L^2(I; X) &\Big|&
v|_{I_n} \in \P_r(I_n; X)\,,\,\,
n=1,\dots,N
\Big\}\,,\quad r \ge 0\,,
\end{array}
\end{equation}
and the continuous in time function space
\begin{equation}
\label{eq:4:3:cGq_Space}
\begin{array}{r@{\,}c@{\,}l@{\,}c@{\,}l}
\mathcal V^{\textnormal{cG($q$)}}_\tau(X) &:=&
\Big\{v \in C(\widebar I; X) &\Big|&
v|_{I_n} \in \P_q(\widebar I_n; X)\,,\,\,
n=1,\dots,N
\Big\}\,,\quad q \ge 1\,.
\end{array}
\end{equation}
which are both discrete and finite dimensional in time.
Functions $v_{\tau}$ of the discrete space
$\mathcal V^{\textnormal{dG($r$)}}_\tau$, or
$\mathcal V^{\textnormal{cG($q$)}}_\tau$,
can be separated into their temporal and spatial parts as
$v_{\tau}(\boldsymbol x, t) = \widehat{v}(\boldsymbol x) \cdot \widetilde{v}(t)$.

For the discretisation in space, we give the partition $\mathcal T_h$
of the domain $\Omega$ into some finite number of disjoint elements $K$,
such that $\widebar \Omega = \bigcup_{K \in \mathcal T_h} \widebar K$.
For simplicity, we choose the elements $K \in \mathcal T_h$ to be
quadrilaterals ($d=2$) or hexahedrals ($d=3$).
We assume that adjacent elements completely share their common edge or face.
We denote by $h_K$ the diameter of the element $K$. The global
space discretisation parameter $h$ is given by $h = \max_{K \in \mathcal T_h} h_K$.
We allow the mesh to be anisotropic but non-degenerated.
For the approximation of the function space
  $\boldsymbol V_0 := \{ \boldsymbol v \in \boldsymbol H(\operatorname{div}; \Omega)
\,|\, (\boldsymbol v \cdot \boldsymbol n)|_{\Gamma_{\fluxC}} = 0 \}$
we use the $\textnormal{RT}_{p-1}$ Raviart-Thomas-N\'ed\'elec finite element space
such that the pressure $\pressure_h$
  is piecewise polynomial of order $(p-1) \ge 0$ on the reference cell and
we use a continuous finite element space $\boldsymbol Q_h^p$,
employing constraints on the degrees of freedom on homogeneous Dirichlet type
boundary on $\Gamma_{\displacement} \times I$,
for vector-valued functions on quadrilaterals or hexahedrals
such that the displacement $\displacement_h$
  is piecewise polynomial of order $p \ge 1$ on the reference cell;
cf. \cite{Koecher2015}.

Furthermore we use the notation of
$\zeta_{n,\jj}(t) \in \P_r(\widebar I_n; \R)$, with $r \ge 0$,
for discontinuous in time basis functions and
$\xi_{n,\jj}(t) \in \P_q(\widebar I_n; \R)$, with $q \ge 1$,
for continuous in time basis functions.
For $p \ge 1$
the vector-valued trial basis functions for the displacement $\displacement_h$
are denoted as $\boldsymbol \chi^K_{\JJ}(\boldsymbol x) \in \boldsymbol Q_h^p$
and the test functions for the mechanics by
$\boldsymbol \omega_h(\boldsymbol x) \in \boldsymbol Q_h^p$,
the vector-valued trial basis functions for the flux $\flux_h$
are denoted as $\boldsymbol \phi^K_{\II}(\boldsymbol x) \in \textnormal{RT}_{p-1}$
and the test functions for Darcy's law by
$\boldsymbol \nu_h(\boldsymbol x) \in \textnormal{RT}_{p-1}$,
the scalar-valued trial basis functions for the pressure $\pressure_h$
are denoted as $\psi^K_{\JJ}(\boldsymbol x) \in \textnormal{RT}_{p-1}$
and the test functions for the evolution equation by
$\mu_h(\boldsymbol x) \in \textnormal{RT}_{p-1}$.
We refer to \cite{Koecher2015} for details on the generation of the basis functions.

\subsection{Discontinuous Galerkin in time method dG($r$)}
\label{sec:4:1:dG}

In this section we present details of the discretisation in time with
a discontinuous Galerkin dG($r$) method with piecewise polynomials
of order $r \ge 0$ in time.
We use here the same spaces
$\mathcal V^{\textnormal{dG($r$)}}_\tau$ from Eq. \eqref{eq:4:2:dGr_Space}
for trial and test functions.
The members of the family of discontinuous Galerkin dG($r$) in time methods are
L-stable, meaning that they are A-stable and, additionally,
the limit of the stability function $\Phi(z)$ of the discretisation of the test
equation $y^\prime=\lambda y$, with $z$ as the product of the
space and time mesh sizes $z=h \tau$, of the respective method is
$\lim_{z \to \infty} \Phi(z) = 0$.
Discontinuous Galerkin methods for the time discretisation usually enforce initial
conditions only weakly, which helps to get a stable solution in the case
of incompatible initial conditions such as jumps of the initial state at $t_0$.
A promising feature is that the generally unknown flux in $t_0$ is not needed at
all and thus the invalidity of Darcy's law from \eqref{eq:2:5:DarcyLaw}
in $t_0^-$ can be neglected.
The fully discrete space-time trial functions are represented by
\begin{equation}
\label{eq:4:4:dGr:fully_discrete_representations}
\begin{array}{rcl@{\quad}rcl}
\displacement_{\tau,h}^{\textnormal{dG($r$)},\s}|_{K \times I_n}(\boldsymbol x,t) &=&
\displaystyle \sum_{\jj=0}^r
  \displacement_h^{\s,n,\jj}(\boldsymbol x)\, \zeta_{n,\jj}(t)\,, &
\displacement_h^{n,\jj}|_K(\boldsymbol x) &=& \displaystyle \sum_{\JJ=1}^{\Nuloc}
  \displacementC_{\s,\JJ,\jj}^{K,I_n}\, \boldsymbol \chi_{\JJ}^K(\boldsymbol x)\,,\\[3.0ex]
\flux_{\tau,h}^{\textnormal{dG($r$)},\s}|_{K \times I_n}(\boldsymbol x,t) &=&
\displaystyle \sum_{\jj=0}^r
  \flux_h^{\s,n,\jj}(\boldsymbol x)\, \zeta_{n,\jj}(t)\,, &
\flux_h^{n,\jj}|_K(\boldsymbol x) &=& \displaystyle \sum_{\JJ=1}^{\Nqloc}
  \fluxC_{\s,\JJ,\jj}^{K,I_n}\, \boldsymbol \phi_{\JJ}^K(\boldsymbol x)\,,\\[3.0ex]
\pressure_{\tau,h}^{\textnormal{dG($r$)},\s}|_{K \times I_n}(\boldsymbol x,t) &=&
\displaystyle \sum_{\jj=0}^r
  \pressure_h^{\s,n,\jj}(\boldsymbol x)\, \zeta_{n,\jj}(t)\,, &
\pressure_h^{n,\jj}|_K(\boldsymbol x) &=& \displaystyle \sum_{\JJ=1}^{\Nploc}
  \pressure_{\s,\JJ,\jj}^{K,I_n}\, \psi_{\JJ}^K(\boldsymbol x)\,.
\end{array}
\end{equation}

The basis functions $\zeta_{n,\jj}(t)$, $\jj=0,\dots,r$, are defined as
Lagrange polynomials via the $(r+1)$-point Gau\ss{} quadrature points
on the reference interval $\widehat I$ and mapped with Eq. \eqref{eq:4:1}.
The assemblies $\alpha^n_{\ii,\jj}$ and $\beta^n_{\ii,\jj}$,
for all $\ii=0,\dots,r$, $\jj=0,\dots,r$ and $n=1,\dots,N$,
for the time discretisation are
\begin{equation}
\label{eq:4:5:dGr:time_assemblies}
\begin{array}{c}
\alpha^n_{\ii,\jj} = 
  \displaystyle \int_{I_n} \zeta_{n,\ii}(t)\, \zeta_{n,\jj}^\prime(t)\, \d t
= \displaystyle \int_{\widehat I} \widehat \zeta_{\ii}(\widehat t)\,
  \widehat \zeta_{\jj}^\prime(\widehat t)\, \d \widehat t\,,
\quad \beta^n_{\ii,\jj} = 
  \displaystyle \int_{I_n} \zeta_{n,\ii}(t)\, \zeta_{n,\jj}(t)\, \d t
= \displaystyle \int_{\widehat I} \widehat \zeta_{\ii}(\widehat t)\,
  \widehat \zeta_{\jj}(\widehat t)\, \tau_n\, \d \widehat t\,,\\[3ex]
\gamma^{n,-}_{\ii} = 
\zeta_{n,\ii}(t_{n-1}^+)
= \widehat \zeta_{\ii}(0)\,,
\quad \gamma^{n,+}_{\ii,\jj} = 
\zeta_{n,\ii}(t_{n-1})\, \zeta_{n,\jj}(t_{n-1})
= \widehat \zeta_{\ii}(0)\, \widehat \zeta_{\jj}(0)\,.
\end{array}
\end{equation}

For the discontinuous Galerkin dG($r$) method in time we replace the
evolution equation of the continuous problem given by the second line of
Eq. \eqref{eq:3:1:FixedStress:Flow} with the following
\begin{equation}
\label{eq:4:6:dGr:jumps_in_evolution_eq}
\begin{array}{@{\!\!\!\!}l@{\;}l}
- & \displaystyle \int_I (w, \divergence \flux^{\s+1}) \,\d t
- \displaystyle \sum_{n=1}^N \int_{I_n} \bigg(w, \bigg( \frac{1}{\BiotModulus}
  + \frac{\BiotCoefficient^2}{\Kdr^*} \bigg) \partial_t \pressure^{\s+1} \bigg) \,\d t
- \displaystyle \bigg(w(t_{n-1}^+), \bigg( \frac{1}{\BiotModulus}
  + \frac{\BiotCoefficient^2}{\Kdr^*} \bigg) [ \pressure^{\s+1} ]_{n-1} \bigg) \\[3.0ex]
= &
\displaystyle
- \int_I (w, f) \,\d t 
+ \int_I (w, \BiotCoefficient\, \strainV(\partial_t \displacement^{\s})) \,\d t
- \displaystyle \sum_{n=1}^N \int_{I_n}
  \bigg(w, \frac{\BiotCoefficient^2}{\Kdr^*}\, \partial_t \pressure^{\s} \bigg) \,\d t
- \bigg(w(t_{n-1}^+), \frac{\BiotCoefficient^2}{\Kdr^*}\, [ \pressure^{\s} ]_{n-1} \bigg)
\end{array}
\end{equation}
and $[ \, \cdot \, ]_n$ denotes a trace operator for the jump
of a discontinuous function in $t_n$ as given by
\begin{equation}
\label{eq:4:7:definition_jump}
\begin{array}{l@{\quad}l@{\quad}l}
\displaystyle [p_{\tau}]_n := p_n^+ - p_n^-\,, &
\displaystyle p_n^- := \lim_{t \to t_n-0} p_{\tau}(t)\,, &
\displaystyle p_n^+ := \lim_{t \to t_n+0} p_{\tau}(t)\,.\\[1.5ex]
\end{array}
\end{equation}

\subsection{Continuous Petrov-Galerkin in time method cG($q$)}
\label{sec:4:2:cG}

In this section we present details of the discretisation in time with
a continuous Petrov-Galerkin cG($q$) method with piecewise polynomials
of order $q \ge 1$ in time.
We use here
$\mathcal V^{\textnormal{cG($q$)}}_\tau$ from Eq. \eqref{eq:4:3:cGq_Space}
for trial functions and
$\mathcal V^{\textnormal{dG($r$)}}_\tau$, $r:=q-1$, from Eq. \eqref{eq:4:2:dGr_Space}
for test functions.
The members of the family of continuous Petrov-Galerkin cG($q$) in time methods
are A-stable but not L-stable.
Continuous Galerkin methods for the time discretisation usually enforce initial
conditions strongly and all solution variables must be continuous on the closure
$\widebar I = [0,T]$.
The fully discrete space-time trial functions are represented by
\begin{equation}
\label{eq:4:8:cGq:fully_discrete_representations}
\begin{array}{rcl@{\quad}rcl}
\displacement_{\tau,h}^{\textnormal{cG($q$)},\s}|_{K \times I_n}(\boldsymbol x,t) &=&
\displaystyle \sum_{\jj=0}^q
  \displacement_h^{\s,n,\jj}(\boldsymbol x)\, \xi_{n,\jj}(t)\,, &
\displacement_h^{n,\jj}|_K(\boldsymbol x) &=& \displaystyle \sum_{\JJ=1}^{\Nuloc}
  \displacementC_{\s,\JJ,\jj}^{K,I_n}\, \boldsymbol \chi_{\JJ}^K(\boldsymbol x)\,,\\[3.0ex]
\flux_{\tau,h}^{\textnormal{cG($q$)},\s}|_{K \times I_n}(\boldsymbol x,t) &=&
\displaystyle \sum_{\jj=0}^q
  \flux_h^{\s,n,\jj}(\boldsymbol x)\, \xi_{n,\jj}(t)\,, &
\flux_h^{n,\jj}|_K(\boldsymbol x) &=& \displaystyle \sum_{\JJ=1}^{\Nqloc}
  \fluxC_{\s,\JJ,\jj}^{K,I_n}\, \boldsymbol \phi_{\JJ}^K(\boldsymbol x)\,,\\[3.0ex]
\pressure_{\tau,h}^{\textnormal{cG($q$)},\s}|_{K \times I_n}(\boldsymbol x,t) &=&
\displaystyle \sum_{\jj=0}^q
  \pressure_h^{\s,n,\jj}(\boldsymbol x)\, \xi_{n,\jj}(t)\,, &
\pressure_h^{n,\jj}|_K(\boldsymbol x) &=& \displaystyle \sum_{\JJ=1}^{\Nploc}
  \pressure_{\s,\JJ,\jj}^{K,I_n}\, \psi_{\JJ}^K(\boldsymbol x)\,.
\end{array}
\end{equation}

The trial basis functions $\xi_{n,\jj}(t)$, $\jj=0,\dots,r,q$, are defined as
Lagrange polynomials via the $(r+1)$-point Gau\ss{} quadrature points
$\widehat t_1, \dots, \widehat t_q$ enhanced by $\widehat t_0=0$
on the reference interval $\widehat I$ and
the test basis functions $\zeta_{n,\ii}(t)$, $\ii=0,\dots,r$, are defined as
Lagrange polynomials via the $(r+1)$-point Gau\ss{} quadrature points
on the reference interval $\widehat I$
and mapped with Eq. \eqref{eq:4:1}.
The assemblies $\alpha^n_{\ii,\jj}$ and $\beta^n_{\ii,\jj}$,
for all $\ii=0,\dots,r$, $\jj=0,\dots,r,q$ and $n=1,\dots,N$,
for the time discretisation are
\begin{equation}
\label{eq:4:9:cGq:time_assemblies}
\begin{array}{c}
\alpha^n_{\ii,\jj} =
  \displaystyle \int_{I_n} \zeta_{n,\ii}(t)\, \xi_{n,\jj}^\prime(t)\, \d t
= \displaystyle \int_{\widehat I} \widehat \zeta_{\ii}(\widehat t)\,
  \widehat \xi_{\jj}^\prime(\widehat t)\, \d \widehat t\,,\quad
\beta^n_{\ii,\jj} =
  \displaystyle \int_{I_n} \zeta_{n,\ii}(t)\,\xi_{n,\jj}(t)\, \d t
= \displaystyle \int_{\widehat I} \widehat \zeta_{\ii}(\widehat t)\,
  \widehat \xi_{\jj}(\widehat t)\, \tau_n\, \d \widehat t\,.
\end{array}
\end{equation}

\subsection{Assemblies from the discretisation in space}
\label{sec:4:3:spatial_assemblies}

We summarise in this subsection the furthermore needed spatial assemblies
of the specific bilinear forms,
that are precisely the $L^2$-inner product integrals over $\Omega$,
of the weak form given by Eqs. \eqref{eq:3:1:FixedStress:Flow},
\eqref{eq:4:6:dGr:jumps_in_evolution_eq}
and
\eqref{eq:3:2:FixedStress:Mech}.
We consider to compile the sparse block matrices
\begin{equation}
\label{eq:4:10:blockmatricesLM}
\boldsymbol L = \begin{pmatrix}
\boldsymbol A & \boldsymbol 0 & \boldsymbol E\\[1ex]
\boldsymbol 0 & \boldsymbol 0 & \boldsymbol B\\[1ex]
\boldsymbol E^\transpose & \boldsymbol B^\transpose & \boldsymbol 0
\end{pmatrix}\,,\quad
\boldsymbol M = \begin{pmatrix}
\boldsymbol 0 & \boldsymbol 0 & \boldsymbol 0\\[1ex]
\boldsymbol 0 & \boldsymbol M_{\flux} & \boldsymbol 0\\[1ex]
\boldsymbol 0 & \boldsymbol 0 & \boldsymbol M_{\pressure}
\end{pmatrix}\,,\quad
\boldsymbol M^{\textnormal{inv}} = \begin{pmatrix}
\boldsymbol 0 & \boldsymbol 0 & \boldsymbol 0\\[1ex]
\boldsymbol 0 & \boldsymbol 0 & \boldsymbol 0\\[1ex]
\boldsymbol 0 & \boldsymbol 0 & \boldsymbol M_{\pressure}^{-1}
\end{pmatrix}\,,
\end{equation}
from the essentially three arising linear equations given by the system
\eqref{eq:2:8:WeakForm}
and needed combinations of spatial trial and test basis functions as well as
occurring derivatives in their fully discrete form
as introduced in Sec. \ref{sec:4:VarSpaceTimeDisc} and
\ref{sec:4:1:dG} for discontinuous in time
or \ref{sec:4:2:cG} for continuous in time approximations.
We assemble the inner block matrices
without the coefficients from a specific time discretisation
such that we can reuse them for different time discretisations,
for fully-coupled monolithic schemes and
for the presented fixed-stress iterative scheme.
We denote by $\Nu$, $\Nq$ and $\Np$ the degrees of freedom in space
for a single temporal degree of freedom and fixed-stress iteration number $\s$
of the fully discrete solutions
$\displacement_{\tau,h}^{\ast,\s}$, $\flux_{\tau,h}^{\ast,\s}$ and
$\pressure_{\tau,h}^{\ast,\s}$, respectively,
with $\ast \in \{ \textnormal{dG}, \textnormal{cG} \}$,
as defined by Sec. \ref{sec:4:1:dG} and \ref{sec:4:2:cG}.
We let the corresponding fully discrete spatial basis functions
have the global numbering from $1$ to $\Nu$, $\Nq$ or, respectively, $\Np$;
cf. \cite{Koecher2015} for details on their generation.
The mechanical element stiffness matrix $\boldsymbol A$ for the corresponding
index sets $1 \le \II \le \Nu$ of test and $1 \le \JJ \le \Nu$ for trial
basis function is given by
\begin{equation}
\label{eq:4:11:matrix:A}
\boldsymbol A = ( a_{\II \JJ} )_{\II \JJ}\,,\quad
a_{\II \JJ} = \displaystyle \sum_{K \in \mathcal{T}_h}
   \displaystyle \sum_{\II=1}^{\Nu}  \displaystyle \sum_{\JJ=1}^{\Nu}
   \displaystyle \int_K
   \strain(\boldsymbol \chi^{\II}) : \elasticity\, \strain(\boldsymbol \chi^{\JJ})
   \,\d \boldsymbol x\,.
\end{equation}
The displacement-pressure coupling matrix $\boldsymbol E$,
using the index sets $1 \le \II \le \Nu$ for test and $1 \le \JJ \le \Np$ for
trial basis functions, is given by
\begin{equation}
\label{eq:4:12:matrix:E}
\boldsymbol E :=
\boldsymbol E_{\displacement \pressure} =
  ( e^{\displacement \pressure}_{\II \JJ} )_{\II \JJ}\,,\quad
e^{\displacement \pressure}_{\II \JJ} = \displaystyle \sum_{K \in \mathcal{T}_h}
   \displaystyle \sum_{\II=1}^{\Nu}  \displaystyle \sum_{\JJ=1}^{\Np}
   \displaystyle \int_K
   \strain(\boldsymbol \chi^{\II}) : \psi^{\JJ} \boldsymbol 1
   \,\d \boldsymbol x\,.
\end{equation}
Since $\strain(\boldsymbol \chi^{\II}) : \boldsymbol 1
= (\nabla \cdot \boldsymbol \chi^{\II})$
we have the pressure-displacement coupling matrix given as
$\boldsymbol E_{\pressure \displacement} =
\boldsymbol E_{\displacement \pressure}^\transpose = \boldsymbol E^\transpose$.
%
%
The matrix $\boldsymbol B$ and
the mass matrices $\boldsymbol M_{\flux}$ and $\boldsymbol M_{\pressure}$
of the Raviart-Thomas-Ned\'el\'ec space are given by
\begin{equation}
\label{eq:4:13:matrix:B}
\boldsymbol B = ( b_{\II \JJ} )_{\II \JJ}\,,\quad
b_{\II \JJ} = \displaystyle \sum_{K \in \mathcal{T}_h}
   \displaystyle \sum_{\II=1}^{\Nq}  \displaystyle \sum_{\JJ=1}^{\Np}
   \displaystyle \int_K
   -(\nabla \cdot \boldsymbol \phi^{\II})\, \psi^{\JJ}
   \,\d \boldsymbol x\,,\quad
\text{with } 1 \le \II \le \Nq \text{ and } 1 \le \JJ \le \Np\,,
\end{equation}
\begin{equation}
\label{eq:4:14:matrix:Mq}
\boldsymbol M_{\flux} = ( m^{\flux}_{\II \JJ} )_{\II \JJ}\,,\quad
m^{\flux}_{\II \JJ} = \displaystyle \sum_{K \in \mathcal{T}_h}
   \displaystyle \sum_{\II=1}^{\Nq}  \displaystyle \sum_{\JJ=1}^{\Nq}
   \displaystyle \int_K
   \boldsymbol \phi^{\II} \cdot \viscosity\, \permeability^{-1}\, \boldsymbol \phi^{\JJ}
   \,\d \boldsymbol x\,,\quad
\text{with } 1 \le \II \le \Nq \text{ and } 1 \le \JJ \le \Nq\,,
\end{equation}
\begin{equation}
\label{eq:4:15:matrix:Mp}
\boldsymbol M_{\pressure} = ( m^{\pressure}_{\II \JJ} )_{\II \JJ}\,,\quad
m^{\pressure}_{\II \JJ} = \displaystyle \sum_{K \in \mathcal{T}_h}
   \displaystyle \sum_{\II=1}^{\Np}  \displaystyle \sum_{\JJ=1}^{\Np}
   \displaystyle \int_K
   \psi^{\II}\, \psi^{\JJ}
   \,\d \boldsymbol x\,,\quad
\text{with } 1 \le \II \le \Np \text{ and } 1 \le \JJ \le \Np\,.
\end{equation}
Remark that $\boldsymbol M_{\pressure}^{-1}$ can be efficiently assembled
since $\boldsymbol M_{\pressure}$ is a block-diagonal or diagonal matrix.
The inverse mass matrix in the pressure space is needed for the system solver
and preconditioning techniques.
%
%
%
The compatibility condition of the initial stress in the mechanical problem
\eqref{eq:3:2:FixedStress:Mech} reads here as
\begin{displaymath}
\displaystyle \sum_{K \in \mathcal{T}_h}
\displaystyle \sum_{\II=1}^{\Nu}
\displaystyle \int_K
  \strain(\boldsymbol \chi^{\II}(\boldsymbol x)) : \stress_0
  \,\d \boldsymbol x
\stackrel{!}{=}
\displaystyle \sum_{K \in \mathcal{T}_h}
\displaystyle \sum_{\II=1}^{\Nu}
\displaystyle \int_K
  \boldsymbol \chi^{\II}(\boldsymbol x) \cdot \density_b\, \gravity(\boldsymbol x, 0)
  \,\d \boldsymbol x
\end{displaymath}
and the projection assembly is given by
$\boldsymbol S_0 = ( s^0_{\II} )_{\II}$, with $1 \le \II \le \Nu$ and
\begin{equation}
\label{eq:4:16:vector:S0}
s^0_{\II} =
\displaystyle \sum_{K \in \mathcal{T}_h}
  \displaystyle \sum_{\II=1}^{\Nu}
  \displaystyle \int_K
  \boldsymbol \chi^{\II}(\boldsymbol x) \cdot \density_b\, \gravity(\boldsymbol x, 0)
  \,\d \boldsymbol x\,.
\end{equation}
%
%
%
The mechanical normal traction assembly is given by
\begin{equation}
\label{eq:4:17:vector:TN}
\boldsymbol T_N^{n,\jj} = ( t^N_{\II} )_{\II}\,,\quad
t^N_{\II} = \displaystyle \sum_{F \in \mathcal{F}_h \cap \Gamma_{\traction}}
   \displaystyle \sum_{\II=1}^{\Nu}
   \displaystyle \int_F
   \boldsymbol \chi^{\II}(\boldsymbol x) \cdot \traction_N(\boldsymbol x, t_{n,\jj})
   \,\d o\,, \quad
\text{with } 1 \le \II \le \Nu\,.
\end{equation}
The gravity assembly for the mechanical problem is given by
\begin{equation}
\label{eq:4:18:vector:Gb}
\boldsymbol G_b^{n,\jj} = ( g^b_{\II} )_{\II}\,,\quad
g^b_{\II} = \displaystyle \sum_{K \in \mathcal{T}_h}
   \displaystyle \sum_{\II=1}^{\Nu}
   \displaystyle \int_K
   \boldsymbol \chi^{\II} \cdot \density_b\, \gravity(\boldsymbol x, t_{n,\jj})
   \,\d \boldsymbol x\,, \quad
\text{with } 1 \le \II \le \Nu\,,
\end{equation}
and the gravity assembly for the flow problem is given by
\begin{equation}
\label{eq:4:19:vector:Gf}
\boldsymbol G_f^{n,\jj} = ( g^f_{\II} )_{\II}\,,\quad
g^f_{\II} = \displaystyle \sum_{K \in \mathcal{T}_h}
   \displaystyle \sum_{\II=1}^{\Nq}
   \displaystyle \int_K
   \boldsymbol \phi^{\II} \cdot \density_f\, \gravity(\boldsymbol x, t_{n,\jj})
   \,\d \boldsymbol x\,, \quad
\text{with } 1 \le \II \le \Nq\,.
\end{equation}
The inhomogeneous pressure boundary assembly for the flow problem is given by
\begin{equation}
\label{eq:4:20:vector:PD}
\boldsymbol P_D^{n,\jj} = ( p^D_{\II} )_{\II}\,,\quad
p^D_{\II} = \displaystyle \sum_{F \in \mathcal{F}_h \cap \Gamma_{\pressure}}
   \displaystyle \sum_{\II=1}^{\Nq}
   \displaystyle \int_F
   (\boldsymbol \phi^{\II}(\boldsymbol x) \cdot \boldsymbol n)\,
   \pressure_D(\boldsymbol x, t_{n,\jj})
   \,\d \boldsymbol x\,, \quad
\text{with } 1 \le \II \le \Nq\,.
\end{equation}
The volumetric source assembly for the flow problem is given by
\begin{equation}
\label{eq:4:21:vector:F}
\boldsymbol F^{n,\jj} = ( f_{\II} )_{\II}\,,\quad
f_{\II} = \displaystyle \sum_{K \in \mathcal{T}_h}
   \displaystyle \sum_{\II=1}^{\Np}
   \displaystyle \int_K
   \psi^{\II}(\boldsymbol x)\, f(\boldsymbol x, t_{n,\jj})
   \,\d \boldsymbol x\,, \quad
\text{with } 1 \le \II \le \Np\,.
\end{equation}
\textbf{Remark.} Within all presented assemblies of the Eqs.
\eqref{eq:4:11:matrix:A} to \eqref{eq:4:21:vector:F}
the innermost loops over the global spatial basis functions
massively reduce in their final implementation due to the compact local support
of the respective finite element basis functions;
cf. for details \cite{Koecher2015}.

\vskip1ex
\textbf{Remark.} The integrals of the Eqs.
\eqref{eq:4:11:matrix:A} to \eqref{eq:4:21:vector:F}
are implemented in a standard way by mapping them on the reference cell
$\widehat K = (0,1)^d$ or reference face $\widehat F = (0,1)^{d-1}$ and
employing a numerical quadrature of necessary order fitted to the occurring
piecewise polynomial order on the reference cell.

\subsection{Fully discrete dG($r$) fixed-stress split}
\label{sec:4:4:fully:dG}

In the following we let $r \ge 0$. Remark that the used matrices, vectors and
coefficients are given by Sec. \ref{sec:4:VarSpaceTimeDisc}, \ref{sec:4:1:dG} and
\ref{sec:4:3:spatial_assemblies}.
The fully discrete dG($r$) in time fixed-stress split reads as:
For $\s = 0,\,\dots,$ do
\begin{enumerate}
\itemsep1ex

\item[(H):] {Find the coefficient vectors $\flux^{\s+1,n,\jj} \in \R^{\Nq}$ and
$\pressure^{\s+1,n,\jj} \in \R^{\Np}$, $0 \le \jj \le r$, from
\begin{equation}
\label{eq:4:22:FixedStress:Flow:dG}
\begin{array}{l@{\!\!\!\! \!\!\!\! \!\!\!\! \!\!\!\!}lcl}
& \displaystyle \sum_{\jj=0}^r \bigg\{
\beta^n_{\ii,\jj} \boldsymbol M_{\flux}\, \flux^{\s+1,n,\jj}
+ \beta^n_{\ii,\jj} \boldsymbol B\, \pressure^{\s+1,n,\jj} \bigg\}
= \displaystyle \sum_{\jj=0}^r \bigg\{
\beta^n_{\ii,\jj} \boldsymbol{G}_f^{n,\jj}
- \beta^n_{\ii,\jj} \boldsymbol{P}_D^{n,\jj}
\bigg\} \,,\\[4.5ex]
& \displaystyle \sum_{\jj=0}^r \bigg\{
\beta^n_{\ii,\jj} \boldsymbol B^\transpose\, \flux^{\s+1,n,\jj}
- (\alpha^n_{\ii,\jj} + \gamma^{n,+}_{\ii,\jj})\,
  \rho_{\pressure}\, \boldsymbol M_{\pressure}\, \pressure^{\s+1,n,\jj}
\bigg\}
=
- \gamma^{n,-}_{\ii}\, \bigg(
  \BiotCoefficient\, \boldsymbol E^\transpose\, \displacement_{n-1}^{-}
  + \frac{1}{\BiotModulus}\, \boldsymbol M_{\pressure}\, \pressure_{n-1}^{-}
\bigg) \\[3ex]
& \quad\quad \displaystyle \sum_{\jj=0}^r \bigg\{
- \beta^n_{\ii,\jj} \boldsymbol F^{n,\jj}
+ (\alpha^n_{\ii,\jj} + \gamma^{n,+}_{\ii,\jj})\, \bigg(
  \BiotCoefficient\, \boldsymbol E^\transpose\, \displacement^{\s,n,\jj}
  - \frac{\BiotCoefficient^2}{\Kdr^*}\, \boldsymbol M_{\pressure}\, \pressure^{\s,n,\jj}
\bigg) \bigg\} \,,
\end{array}
\end{equation}
for $0 \le \ii \le r$ and
where we denote by $\displacement_{n-1}^{-} \in \R^{\Nu}$ and by
$\pressure_{n-1}^{-} \in \R^{\Np}$ the vectors of the interpolation or
projection of the respective initial value function or solution at $t_{n-1}^-$.
}

\item[(M):] {%
Find the coefficient vectors $\displacement^{\s+1,n,\jj} \in \R^{\Nu}$,
$0 \le \jj \le r$, from
\begin{equation}
\label{eq:4:23:FixedStress:Mech:dG}
\begin{array}{l}
\displaystyle \sum_{\jj=0}^r \bigg\{
\beta^n_{\ii,\jj}\, \boldsymbol A\,
  (\displacement^{\s+1,n,\jj} - \boldsymbol U_0) \bigg\} = \\[3ex] 
\quad\quad\quad
\displaystyle \sum_{\jj=0}^r \bigg\{
  \beta^n_{\ii,\jj}\, \boldsymbol G_b^{n,\jj}
+ \beta^n_{\ii,\jj}\, \boldsymbol T_N^{n,\jj}
- \beta^n_{\ii,\jj}\, \boldsymbol S_0
\bigg\} + 
%
\displaystyle \sum_{\jj=0}^r \bigg\{
  \beta^n_{\ii,\jj}\, \BiotCoefficient\, \boldsymbol E\,
    (\pressure^{\s+1,n,\jj} - \boldsymbol P_0)
\bigg\} \,,
\end{array}
\end{equation}
for $0 \le \ii \le r$ and
where we denote by $\boldsymbol U_0 \in \R^{\Nu}$ and by
$\boldsymbol P_0 \in \R^{\Np}$ the vectors of the interpolation or
projection of the respective initial value function at $t_0$.
}
\end{enumerate}
until convergence and then marching forwardly through all time subintervals $I_n$.

\vskip1ex
\textbf{Remark.} For the construction of appropriate iterative linear solvers
and preconditioning techniques as well as their implementation we refer to
\cite{Koecher2018z,Koecher2018y,Koecher2016a,Koecher2016,Koecher2015}.

\subsection{Fully discrete cG($q$) fixed-stress split}
\label{sec:4:5:fully:cG}

In the following we let $q \ge 1$ and put $r=q-1$.
Remark that the used matrices, vectors and coefficients are given by
Sec. \ref{sec:4:VarSpaceTimeDisc}, \ref{sec:4:2:cG} and
\ref{sec:4:3:spatial_assemblies}.
The fully discrete cG($q$) in time fixed-stress splitting scheme reads as:
For $\s = 0,\,\dots,$ do
\begin{enumerate}
\itemsep1ex

\item[(H):] {Find the coefficient vectors $\flux^{\s+1,n,\jj} \in \R^{\Nq}$ and
$\pressure^{\s+1,n,\jj} \in \R^{\Np}$, $1 \le \jj \le q$, from
\begin{equation}
\label{eq:4:24:FixedStress:Flow:cG}
\begin{array}{lcl}
\displaystyle \sum_{\jj=0}^q \bigg\{
\beta^n_{\ii,\jj} \boldsymbol M_{\flux}\, \flux^{\s+1,n,\jj}
+ \beta^n_{\ii,\jj} \boldsymbol B\, \pressure^{\s+1,n,\jj} \bigg\}
= \displaystyle \sum_{\jj=0}^q \bigg\{
\beta^n_{\ii,\jj} \boldsymbol{G}_f^{n,\jj}
- \beta^n_{\ii,\jj} \boldsymbol{P}_D^{n,\jj}
\bigg\} \,,\\[4.5ex]
\displaystyle \sum_{\jj=0}^q \bigg\{
\beta^n_{\ii,\jj} \boldsymbol B^\transpose\, \flux^{\s+1,n,\jj}
- \alpha^n_{\ii,\jj}\,
  \rho_{\pressure}\, \boldsymbol M_{\pressure}\, \pressure^{\s+1,n,\jj}
\bigg\}
= \\[3ex]
\quad\quad \displaystyle \sum_{\jj=0}^q \bigg\{
- \beta^n_{\ii,\jj} \boldsymbol F^{n,\jj}
+ \alpha^n_{\ii,\jj}\, \bigg(
  \BiotCoefficient\, \boldsymbol E^\transpose\, \displacement^{\s,n,\jj}
  - \frac{\BiotCoefficient^2}{\Kdr^*}\, \boldsymbol M_{\pressure}\, \pressure^{\s,n,\jj}
\bigg) \bigg\} \,,
\end{array}
\end{equation}
for $0 \le \ii \le r$ and
where we denote
by $\displacement^{\star,n,0} \in \R^{\Nu}$,
by $\flux^{\star,n,0} \in \R^{\Nq}$ and
by $\pressure^{\star,n,0} \in \R^{\Np}$ the vectors of the interpolation or
projection of the respective initial value function or solution at $t_{n-1}$.
}

\item[(M):] {%
Find the coefficient vectors $\displacement^{\s+1,n,\jj} \in \R^{\Nu}$,
$1 \le \jj \le q$, from
\begin{equation}
\label{eq:4:25:FixedStress:Mech:cG}
\begin{array}{l}
\displaystyle \sum_{\jj=0}^q \bigg\{
\beta^n_{\ii,\jj}\, \boldsymbol A\,
  (\displacement^{\s+1,n,\jj} - \boldsymbol U_0) \bigg\} = \\[3ex] 
\quad\quad\quad
\displaystyle \sum_{\jj=0}^q \bigg\{
  \beta^n_{\ii,\jj}\, \boldsymbol G_b^{n,\jj}
+ \beta^n_{\ii,\jj}\, \boldsymbol T_N^{n,\jj}
- \beta^n_{\ii,\jj}\, \boldsymbol S_0
\bigg\} + 
%
\displaystyle \sum_{\jj=0}^q \bigg\{
  \beta^n_{\ii,\jj}\, \BiotCoefficient\, \boldsymbol E\,
    (\pressure^{\s+1,n,\jj} - \boldsymbol P_0)
\bigg\} \,,
\end{array}
\end{equation}
for $0 \le \ii \le r$ and
where we denote by $\boldsymbol U_0 \in \R^{\Nu}$ and by
$\boldsymbol P_0 \in \R^{\Np}$ the vectors of the interpolation or
projection of the respective initial value function at $t_0$ and
where we denote
by $\displacement^{\star,n,0} \in \R^{\Nu}$ and
by $\pressure^{\star,n,0} \in \R^{\Np}$ the vectors of the interpolation or
projection of the respective initial value function or solution at $t_{n-1}$.
}
\end{enumerate}
until convergence and then marching forwardly through all time subintervals $I_n$.

\vskip1ex
\textbf{Remark.} For the construction of appropriate iterative linear solvers
and preconditioning techniques as well as their implementation we refer to
\cite{Koecher2018z,Koecher2018y,Koecher2016a,Koecher2016,Koecher2015}.

\subsection{Fully discrete dG($r$)-cG($q$) fixed-stress split}
\label{sec:4:6:dGcGq}

In the following we introduce a new Galerkin time discretisation for the
fixed-stress operator splitting scheme for the Biot's system of poroelasticity
by mixing
the discontinuous Galerkin dG($r$), $r \ge 0$, discretisation of Sec.
\ref{sec:4:1:dG} and
the continuous Galerkin cG($q$), $q \ge 1$, discretisation of Sec. \ref{sec:4:2:cG}
as time discretisation of the continuous space-time problem
given by Sec. \ref{sec:3:OperatorSplitting:Continuum}
with the Eqs. \eqref{eq:3:1:FixedStress:Flow} - \eqref{eq:3:2:FixedStress:Mech}.
In Sec. \ref{sec:5}, we will show numerically the applicability and
numerical stability of the new method and refer a formal proof of the stability
and convergence to a forthcoming work.
We introduce the set of time intervals
$\mathcal T_\tau := \{ I_n \subset I \,|\, n=1,\dots,N \}$ with notation introduced
in Sec. \ref{sec:4:VarSpaceTimeDisc} for the partition of the global time interval
$I=(0,T)$ into disjoint subintervals $I_n$.
We divide the set $\mathcal T_\tau$ into disjoint subsets
$\mathcal T_\tau^{\text{dG}}$ and $\mathcal T_\tau^{\text{cG}}$
for marking time subintervals $I_n$ in which the solutions will be approximated
with the family of discontinuous and continuous Galerkin time discretisations,
respectively.
The partition of $\mathcal T_\tau$ must hold
$$
\mathcal T_\tau^{\text{dG}} \cup \mathcal T_\tau^{\text{cG}} =
  \mathcal T_\tau = \{ I_1, \dots, I_N \}
\quad\text{and}\quad
\mathcal T_\tau^{\text{dG}} \cap \mathcal T_\tau^{\text{cG}} =
  \emptyset\,.
$$
We remark that the time subintervals in $\mathcal T_\tau^{\text{dG}}$ and
$\mathcal T_\tau^{\text{cG}}$ must not necessarily contain continuous subranges
of time subintervals, e.g. a partition
$\mathcal T_\tau^{\text{dG}} = \{ I_1,\, I_4, \dots, I_7 \}$ and
$\mathcal T_\tau^{\text{cG}} = \{ I_2, I_3,\, I_8, \dots, I_N \}$
is possible for instance.
The new approach still make use of a time marching process,
in which decoupled problems are solved completely on a time subinterval $I_n$
before solving the next problem on $I_{n+1}$.
Within continuous subranges of time subintervals,
  i.e. e.g. $\{ I_{n_1}, I_{n_1+1}, \dots, I_{n_2} \}$,
  with $1 \le n_1 \le \dots \le n_2 \le N$,
  in either $\mathcal T_\tau^{\text{dG}}$ or $\mathcal T_\tau^{\text{cG}}$,
we use the the corresponding fixed-stress scheme as given in Sec. \ref{sec:4:1:dG}
or \ref{sec:4:2:cG}, respectively.
Whenever two different types of time discretisations appear between two
consecutive time subintervals, e.g.
$I_{n_3+1}=I_{n_4}$, with $1 \le \dots \le n_3 < n_3+1 = n_4 \le \dots \le N$,
and $I_{n_3} \in \mathcal T_\tau^{\ast}$
having $\ast \in \{ \text{dG}, \text{cG} \}$ and
$I_{n_4}$ in the respective other subset of $\mathcal T_\tau$,
we interpolate the current final time $t_{n_3}$ solutions as initial values
for the next problem on $I_{n_4}$ in the sense of a time marching scheme.

\vskip1ex
\textbf{Remark.} We formally allow different polynomial degrees and time
subinterval lengths in every $I_n$ in the sense of a corresponding one-dimensional
$hp$-adaptive method known from variational space discretisations
for our new time discretisation.

\vskip1ex
\textbf{Remark.} Mixing continuous and discontinuous Galerkin methods for
the spatial discretisation of several problems has already been done;
cf. e.g. \cite{Dawson2002} and references therein.

\vskip1ex
\textbf{Scheme 1 (dG(1)-cG(1) for incompatible initial data).}
In Sec. \ref{sec:5}, we use of the partition
  $\mathcal T_\tau^{\text{dG($1$)}} = \{ I_1 \}$ and
  $\mathcal T_\tau^{\text{cG($1$)}} = \{ I_2, \dots, I_N \}$
using a piecewise linear approximation in time for both time discretisation
families.

\section{Numerical Examples}
\label{sec:5}

In the following section we study the stability properties and
computational efficiency of the introduced Galerkin time discretisations
of discontinuous and continuous families of Sec. \ref{sec:4:1:dG} and
Sec. \ref{sec:4:2:cG} as well as the new mixed discontinuous-continuous Galerkin
time discretisation dG(1)-cG(1) for incompatible initial data given by
Scheme 1 of Sec. \ref{sec:4:6:dGcGq}.

\begin{figure}
\centering
\subfloat[Schematic test setting.]{%
\includegraphics[width=0.3\linewidth]{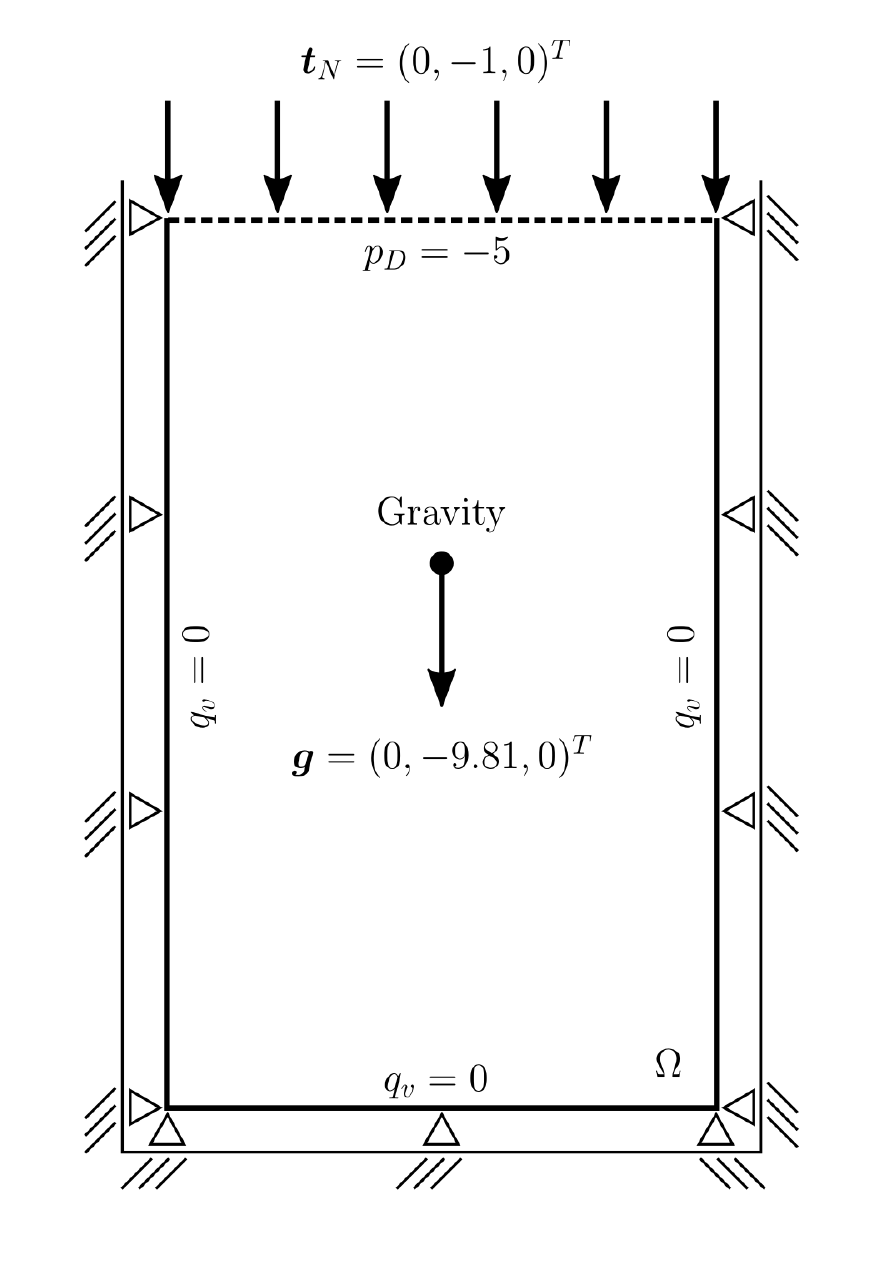}
}
%
\quad\quad
\subfloat[Solution for $t=0.1$.]{%
\includegraphics[width=0.45\linewidth]{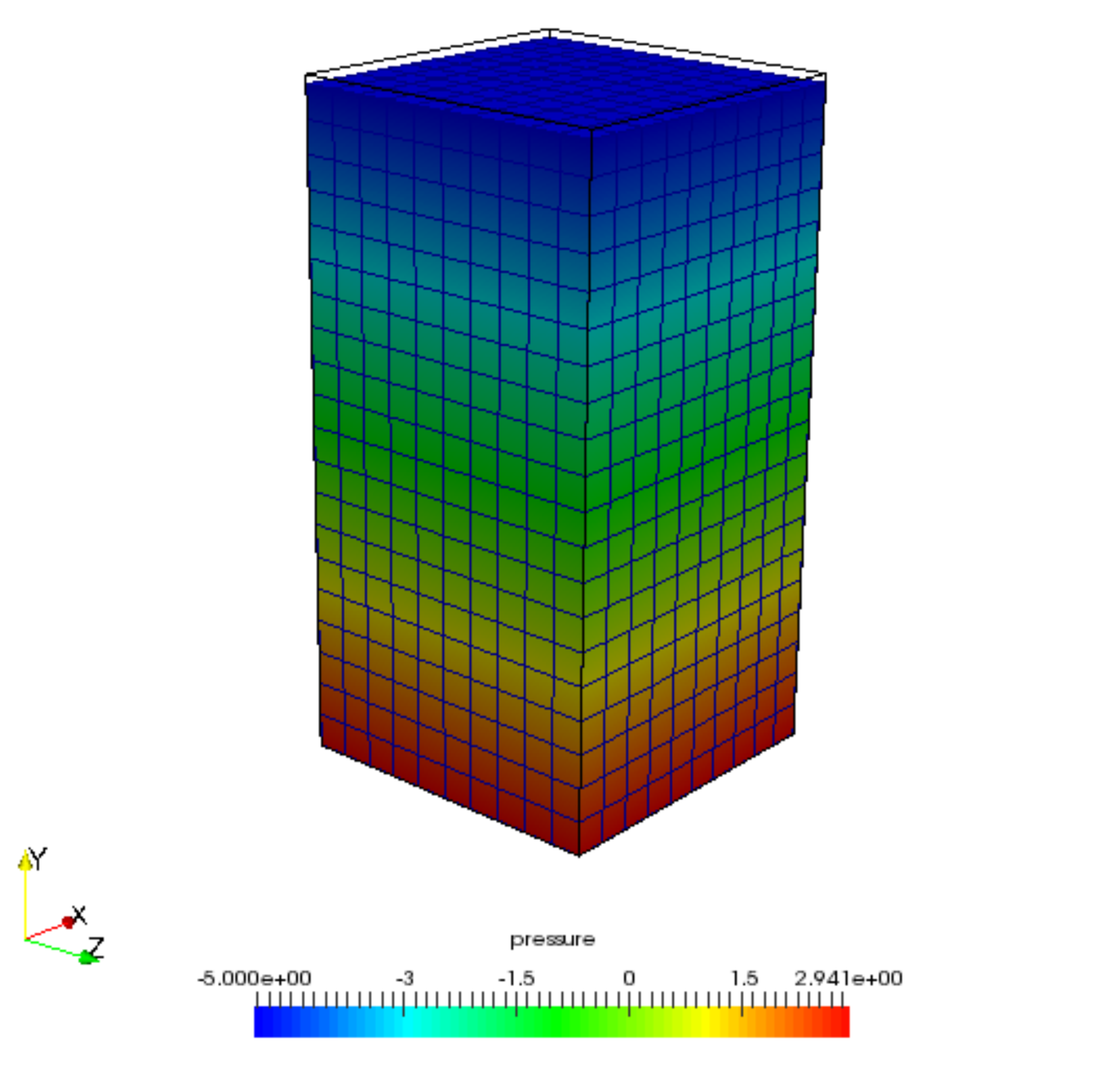}
}
\caption{Schematic test setting given by a two-dimensional cut illustration (a)
and solution illustration (b) for Sec. \ref{sec:5}.
In (b) the displacement field $\displacement$ is illustrated by warping the domain
$\Omega$ and the pressure difference distribution $\pressure$ by the colour
at $t=0.1$ of a three-dimensional Terzaghi-type consolidation test problem
with piecewise quadratic displacement and piecewise linear pressure
approximation in space and the new dG($1$)-cG($1$) fixed-stress solver
for incompatible initial data given by Scheme 1 in Sec. \ref{sec:4:6:dGcGq}.}
\label{fig:1}
\end{figure}

We use here a Terzaghi-type consolidation test enforcing an uniaxial compression
on a three-dimensional fluid-saturated poroelastic block placed in an aquarium
with non-deformable walls and open top surface
which is loaded by gravity and, moreover, instantaneously in $t_0=0^+$
with an outer normal traction and therefore has incompatible initial data.
The schematic test setting of all numerical examples is illustrated
with a two-dimensional cut by Fig. \ref{fig:1} (a) and further explained in detail.

We remark that the overall test setting yields a problem which ensures the
existence and uniqueness of solutions even with the componentwise taken
mixed Dirichlet and Neumann type boundary conditions
for the mechanical problem.
The contact problem of the mechanical problem,
occuring here by hitting the non-deformable walls of the fluid-saturated
poroelastic material under load,
is resolved here, for simplicity and as usual, by setting physically reasonable
homogeneous Dirichlet type boundary data on the displacement field in a
componentwise taken sense.

Remark that the test setting yields benchmark problem,
since it was already used in \cite{Jha2007,Koecher2016,Cryer1963,Terzaghi1943}
in corresponding an one- or two-dimensional way,
but we use here fully three-dimensional setting.
Remark also that for the generation of the solutions of the displacement and the
pressure for the given test problem even the classical one-dimensional test is
enough to get satisfying results of physically relevance.
Anyhow, we use here the corresponding three-dimensional form since this stresses
much more the applied space-time discretisations, the applied fixed-stress
iterative solvers as well as the linear solver and preconditioning techniques
underneath to show advantages and disadvantages of the approaches,
which is in the focus of this work and is of further interest in the current
research,
and to simplify the verification by comparision with known results of the
computed solutions for the reader.

The pressure solution below the constrained top surface includes a
Mandel-Cryer effect,
which is a pressure overshoot shortly after the initial time
since the fluid is not able to escape instantaneously;
cf. \cite{Cryer1963,Mandel1953}.
This Mandel-Cryer effect is an additional difficulty for non-robust solver
strategies as they can appear for iteratively-coupled solvers,
e.g. fixed-stress solver technologies for higher-order time discretisations
employing a complicated nested Schur complement solver technology,
which will be also studied by the following examples.
We note that such difficulties can be overcome by using
monolithic solver strategies instead of iteratively coupled solvers,
but they have the need for an efficient and scalable preconditioning technology
which will be studied in our forthcoming work in \cite{Koecher2018z}.

All numerical results are computed with the \texttt{DTM++/biot-fs} solver suite
which is a frontend solver of the author for the \texttt{deal.II} v8.5.1 library;
cf. \cite{Koecher2018y,Koecher2015,Koecher2016} and \cite{dealII85}.
They enforce 10 distributed parallel (MPI) processes each
for the parallelisation and distribution of the spatial domain such that each
process locally owns 25 mesh cells $K$.
Note that different numbers of overall degrees of freedoms need to be computed
for different time discretisations of the same piecewise polynomial order
in time for the introduced families.
In general, discontinuous Galerkin (time) discretisations are
more expensive in any terms for the same piecewise polynomial degree
as their continuous Galerkin discretisation counterparts.

We optimise all fixed-stress iterative schemes by using $\Kdr^* = \lambda+2\mu$,
with $\lambda$ and $\mu$ being the first and second Lam\'e parameters as introduced
in Sec. \ref{sec:2:QuasiStaticBiotSystem},
as the most optimal $\Kdr^*$ value under uniaxial loading
without deeper knowledge on the problem as it is studied in detail in
\cite{Koecher2017a}.
We note that the numbers of the fixed-stress iterations is not a measurement of
the performance in this work since it has to be shown to lead to non-comparable
results due to different tunings.
All numerical solvers are tuned in a way that they deliver
comparable results in terms of computational efficiency and accuracy.

For the following experiments we let the test setting illustrated by
Fig. \ref{fig:1} (a) being the same for all numerical experiments.
The spatial domain is set to $\Omega = (0,0.5) \times (0,1) \times (0,0.5)$
with $h=\sqrt{0.03}$ and time interval is $I=(0, 0.5)$ with $\tau = \tau_n = 0.001$.
On the top surface, i.e. precisely
$\Gamma_{\pressure} = \{ (x, 1, z)^\transpose \in \R^3 \,|\,
(x,z)^\transpose \in (0,0.5)^2 \}$,
we prescribe an inhomogeneous pressure boundary condition with
$\pressure_D = -5$ and a zero flux boundary condition elsewhere
for the flow problem.
Note that our global pressure solutions and boundary conditions here must be
understand as difference to the corresponding reference state.
For the mechanical problem we prescribe an uniform
compressive traction $\traction_N = (0,-1,0)^\transpose$ on the top surface
and componentwise mixed homogeneous Dirichlet and homogeneous traction conditions
elsewhere.
The piecewise polynomial approximation degree in space is choosen such that we
have a piecewise triquadratic approximation of the displacement field
($\boldsymbol Q_h^2$ finite element space)
and a piecewise trilinear approximation of the pressure
($\text{RT}_1$ mixed finite element space).
For all material parameters we use the values given in \cite[Sec. 5.1]{Jha2007}.

In Fig. \ref{fig:1} (b) the displacement field $\displacement$ is illustrated
by warping the domain $\Omega$ and the pressure difference distribution $\pressure$
by the colour at $t=0.1$ for the new dG($1$)-cG($1$) fixed-stress solver
for incompatible initial data given by Scheme 1 in Sec. \ref{sec:4:6:dGcGq}.
For better comparision of the fixed-stress solvers with different time discretisations
in the following Sec. \ref{sec:5:1} and \ref{sec:5:2} we compute the mean value
of the pressure solution along the embedded one-dimensional line
$\boldsymbol c_L(s)=(s,0.15,0.25)^\transpose$ with $s \in (0,0.5)$. Note that
the line $\boldsymbol c_L$ is aligned to a constant pressure layer.

The further structure of this section reads as follows.
In Sec. \ref{sec:5:1} we compare the time discretisations dG(0), cG(1) and dG(1),
in Sec. \ref{sec:5:2} we study the new dG(1)-cG(1) for incompatible initial data
given by Scheme 1 of Sec. \ref{sec:4:6:dGcGq} against the dG(0) and dG(1)
time discretisations
and in Sec. \ref{sec:5:3} we study the performance of all fixed-stress solvers.

\input{Koecher2018-mixed-splitting-compare-pressure-51}

\input{Koecher2018-mixed-splitting-compare-pressure-52}

\subsection{Test case: dG(0), cG(1) and dG(1) for incompatible initial data}
\label{sec:5:1}

In the first numerical study we compare (a) the dG(0) with the cG(1) time
discretisations,
which are the members of lowest-order type in time in their family,
and (b) the dG(0) and the dG(1) time discretisations.
The corresponding pressure solutions averaged on $\boldsymbol c_L$ are illustrated
in Fig. \ref{fig:2:compare:pressure:dG0_dG1_cG1} (a) and (b).

The dG(0) time discretisation can be re-identified
as the classical backward Euler time discretisation, which is known as robust
time discretisation, but it delivers only piecewise constant approximation in time.
The numerical solution of the dG(0) fixed-stress solver match with the known data
from the literature; cf. \cite{Jha2007,Koecher2016}.

In contrast, the cG(1) fixed-stress solver produces an unphysically oscillating
and therefore unsuitable numerical solution. Note that even a monolithic solver
for the presented cG(1) time discretisation has the same issue with unphysical
oscillations in the case of incompatible data.
This result of the cG(1) solver is unsatisfying since it has a piecewise linear
approximation in time for comparable computational costs as the dG(0) solver
which will be analysed in Sec. \ref{sec:5:3}.

Secondly, in Fig. \ref{fig:2:compare:pressure:dG0_dG1_cG1} (b)
we compare the dG(0) and dG(1) fixed-stress solver solutions.
For $t \in (0, 0.045)$ the pressure results are only slightly different,
but after the Mandel-Cryer effect the dG(1) fixed-stress solver gets into trouble,
which results in unforeseeable number of fixed-stress iterations.
As we use here the dG(1) fixed-stress solver without specific tuning to the
test problem to allow for comparable results for Sec. \ref{sec:5:3},
the pressure solution clearly shows imprecise results for $t > 0.045$.
Note that this issue can be resolved by enforcing higher tolerance values within
the complex and nested linear system solver
but then the simulation wall time will increase massively.
Remark that we have presented the needed complex inner solver
structure for the dG(1) fixed-stress solver in \cite{Koecher2016,Koecher2016a}.

\subsection{Test case: dG(1)-cG(1) against dG(0) and dG(1) for incompatible initial data}
\label{sec:5:2}

In the second numerical study we compare (a) the dG(0) time discretisation with
the new dG(1)-cG(1) coupled time discretisation given by Scheme 1 in Sec.
\ref{sec:4:6:dGcGq} and (b) the dG(1) time discretisation with the new dG(1)-cG(1)
time discretisation.
The corresponding pressure solutions averaged on $\boldsymbol c_L$ are illustrated
in Fig. \ref{fig:3:compare:pressure:dG0_dG1_dG1cG1} (a) and (b).

In contrast to Sec. \ref{sec:5:1},
the new dG(1)-cG(1) fixed-stress solver computes suitable numerical solution
which matches mostly with the dG(0) solution as it is illustrated in Fig.
\ref{fig:3:compare:pressure:dG0_dG1_dG1cG1} (a).
Note that the slight difference in the solutions are expected from different
convergence orders of the time discretisations coming from the piecewise constant
and, respectively, piecewise linear approximations in time.
This hypothesis has to be shown with rigorous error analysis in a future work.
Additionally remark that the spatial discretisation used for all test problems
was initially converged such that only differences for the time discretisation
can be clearly presented here.

In Fig. \ref{fig:3:compare:pressure:dG0_dG1_dG1cG1} (b) we study the piecewise
linear approximations in time with the global dG(1) fixed-stress solver and
the new dG(1)-cG(1) fixed-stress solver.
The pressure solutions of the new dG(1)-cG(1) solver fit to the solutions of
the dG(1) fixed-stress solver for $t < 0.045$.
The new approach does not match with non-tuned dG(1) fixed-stress solver
for $t > 0.045$,
which let us lead to interpret the overall dG(1) numerical solution as unsatisfying
result for the test problem and solver tolerances here.
Remark that we have studied in \cite{Koecher2016} that also the dG(1) fixed-stress
solver gives satisfying numerical approximations.

\subsection{Performance comparision: dG(0), dG(1)-cG(1), cG(1) and dG(1)}
\label{sec:5:3}

In this subsection we analyse the computational performance and efficiency
of the fixed-stress solvers dG(0), cG(1), dG(1) and the new dG(1)-cG(1)
used in Sec. \ref{sec:5:1} and \ref{sec:5:2} by giving the simulation wall times.
Note that not the simulation wall time (sometimes denoted as CPU time)
itself is a measurement of the efficiency, instead the difference between the
wall times of the solvers demonstrate the efficiency.
We remark again that our implementation in the \texttt{DTM++/biot-fs}
solver suite leads to comparable solvers since most of the code is reused
in a hierarchical and modularised way; cf. \cite{Koecher2018y}.

The simulation wall times for 
the fixed-stress solvers dG(0), cG(1), dG(1) and the new dG(1)-cG(1)
fixed-stress solvers differ in a major way:
\begin{itemize}
\item dG(0) fixed-stress solver: 56 seconds,
\item dG(1)-cG(1) fixed-stress solver: 90 seconds,
\item cG(1) fixed-stress solver: 184 seconds and
\item dG(1) fixed-stress solver: 1573 second.
\end{itemize}
to fully compute the problem given in Sec. \ref{sec:5} for $I=(0,0.5)$ with
$\tau = \tau_n = 0.001$ for all solvers and including all preprocessing and
postprocessing steps such as parsing the input parameters, grid generation,
assembly as well as data output of the global solution.
Note that the real cpu time, not the simulation wall time, can be computed
by multiplying the respective wall time by the factor of 10 since we have used
10 MPI processes for each simulation.
Again, as it is studied in detail in Sec. \ref{sec:5:1} and \ref{sec:5:2},
the solution of the cG(1) fixed-stress solver is not feasible to be used
while the solution of the new dG(1)-cG(1) approach given by Scheme 1 in
Sec. \ref{sec:4:6:dGcGq} yields accurate numerical solutions.
Moreover, the new approach performs the complete simulation with mostly the
same simulation wall time as the dG(0) fixed-stress solver. We remark that the
gap between the wall times comes from solving with the dG(1) fixed-stress solver
on $I_1$ in the new approach.

We finally remark that the obvious lack of accuracy, performance and efficiency
of the dG(1)-fs solver comes here due to the fixed-stress iterative approach,
in a forthcoming work we show that this can be resolved by a monolithic solver approach
using an efficient and flexible preconditioning technology;
\cite{Koecher2018z}.

\section{Conclusions}
\label{sec:6}

In this work we studied in detail Biot's quasi-static system of poroelasticity
which we have established from conservation and balance laws and other physical
principles under necessary and natural assumptions on the data.
Then we presented the fixed-stress iterative approach for decoupling the flow
and the mechanical subproblems for the infinite-dimensional space-time problem.
Further we studied the application of discontinuous and continuous Galerkin
time discretisations of arbitrary polynomial order in time.
We presented details for the assemblies of the mixed discretisation in space
in a way such that they can be reused for an arbitrary time discretisation and
for splitting or monolithic schemes.
A new dG($r$)-cG($r$) in time approach, using dG($r$) and cG($r$) time
discretisations in different time regions, was introduced.
With a three-dimensional numerical experiment of physical relevance we studied
the stability and computational efficiency of the presented schemes.
The new dG($1$)-cG($1$) approach has been shown numerically to be as robust and
computationally efficient as the classically used backward Euler time
discretisation while the global (higher-order) dG($1$) in time fixed-stress approach
did not yield satisfying results.
The application of monolithic solvers instead of iteration schemes will be
done in a forthcoming work; \cite{Koecher2018z}.
The stability and convergence of the new dG($1$)-cG($1$) in time approach
needs to be rigorously analysed in a future work.


\end{document}